\documentclass[11pt]{amsart}

\usepackage{amsfonts}
\usepackage{amsmath}
\usepackage{amsthm}
\usepackage{graphicx}
\usepackage{caption}
\usepackage{natbib}
\usepackage{color}
\usepackage{hyperref}
\usepackage{setspace}
\usepackage{natbib}
\usepackage[T1]{fontenc}
\usepackage{times}
\usepackage[left=2cm,right=2cm,top=2.5cm,bottom=2.5cm]{geometry}

\usepackage{bigints}

\newcommand{\EXCLUDE}[1]{}
\newcommand{\no}{\nonumber}

\newcommand{\pr}[1]{\mathbb{P}\left\{ #1 \right\}}

\newcommand{\EXP}[1]{\mathbb{E}\!\left[#1\right] }

\newcommand{\VAR}[1]{\mathsf{VAR}\!\left(#1\right) }
\newcommand{\remove}[1]{}

\newcommand{\beq}{\begin{eqnarray}}
	\newcommand{\eeq}{\end{eqnarray}}
\newcommand{\beqq}{\begin{eqnarray*}}
	\newcommand{\eeqq}{\end{eqnarray*}}
\def\:{:\,}

\newtheorem{theorem}{Theorem}[section]
\newtheorem{corollary}[theorem]{Corollary}
\newtheorem{lemma}[theorem]{Lemma}
\newtheorem{proposition}[theorem]{Proposition}
\newtheorem{remark}[theorem]{\bf Remark}
\newtheorem{definition}[theorem]{\bf Definition}

\def\Cech{\v{C}ech\ }

\def\VR{Vietoris-Rips\ }

\newcommand{\ep}{\epsilon}
\newcommand{\al}{\alpha}

\newcommand{\del}{\delta}

\newcommand{\sg}{\sigma}

\def\th{\theta}

\def\ze{\zeta}

\def\mR{\mathbb{R}}

\def\mZ{\mathbb{Z}}

\def\mN{\mathbb{N}}

\newcommand{\md}{{\,d}}

\newcommand{\cP}{\mathcal{P}}
\newcommand{\cU}{{\mathcal U}}
\newcommand{\cI}{{\mathcal I}}
\newcommand{\cC}{{\mathcal C}}
\newcommand{\cD}{{\mathcal D}}
\newcommand{\cR}{{\mathcal R}}
\newcommand{\cX}{{\mathcal X}}

\newcommand{\cN}{{\mathcal N}}

\newcommand{\cK}{{\mathcal K}}

\def\1{\mathbf{1}}

\def\x{\bold{x}}
\def\y{\bold{y}}
\def\z{\bold{z}}
\def\bv{\bold{v}}
\def\bze{\bold{\ze}}
\def\X{\bold{X}}

\def\hc{h^{\cC}}
\def\hr{h^{\cR}}
\def\Qcu{Q^{\cC,\cU}}
\def\Qcd{Q^{\cC,\cD}}
\def\Qru{Q^{\cR,\cU}}
\def\Qrd{Q^{\cR,\cD}}

\def\Gcu{G^{\cC,\cU}}
\def\Gcd{G^{\cC,\cD}}

\def\Gru{G^{\cR,\cU}}
\def\Grd{G^{\cR,\cD}}

\def\lb{\left(}
\def\rb{\right)}

\def\lab{\label}

\def\nn{\nonumber}
\def\f{\frac}

\def\ls{\displaystyle{\limsup_{n \to \infty}\;}}

\def\lo{\displaystyle{\lim_{n \to \infty}\;}}

\begin{document}
	
	
	\title[Isolated faces in random geometric complexes]{Thresholds for vanishing of `Isolated' faces in random \v{C}ech and Vietoris-Rips complexes.}
\author{Srikanth K. Iyer}
\address{(SKI) Department of Mathematics\\
		Indian Institute of Science\\
		Bangalore, 	India.}
\email{skiyer@iisc.ac.in}
	
\author{D. Yogeshwaran}
\address{(DY) Theoretical Statistics and Mathematics unit\\
	Indian Statistical Institute \\
	Bangalore, India}
\thanks{DY's research was supported in part by DST-INSPIRE Faculty fellowship and CPDA from the Indian Statistical Institute.}
\email{d.yogesh@isibang.ac.in}%
%
\begin{abstract}
We study combinatorial connectivity for two models of random geometric complexes. These two models - \v{C}ech and Vietoris-Rips complexes - are built on a homogeneous Poisson point process of intensity $n$ on a $d$-dimensional torus using balls of radius $r_n$. In the former, the $k$-simplices/faces are formed by subsets of $(k+1)$ Poisson points such that the balls of radius $r_n$ centred at these points have a mutual interesection and in the latter, we require only a pairwise intersection of the balls. Given a (simplicial) complex (i.e., a collection of $k$-simplices for all $k \geq 1$), we can connect $k$-simplices via $(k+1)$-simplices (`up-connectivity') or via $(k-1)$-simplices (`down-connectivity). Our interest is to understand these two combinatorial notions of connectivity for the random \v{C}ech and Vietoris-Rips complexes asymptically as $n \to \infty$. In particular, we analyse in detail the threshold radius for vanishing of isolated $k$-faces for up and down connectivity of both types of random geometric complexes. Though it is expected that the threshold radius $r_n = \Theta((\frac{\log n}{n})^{1/d})$ in coarse scale, our results give tighter bounds on the constants in the logarithmic scale as well as shed light on the possible second-order correction factors. Further, they also reveal interesting differences between the phase transition in the \v{C}ech and Vietoris-Rips cases. The analysis is interesting due to the non-monotonicity of the number of isolated $k$-faces (as a function of the radius) and leads one to consider `monotonic' vanishing of isolated $k$-faces. The latter coincides with the vanishing threshold mentioned above at a coarse scale (i.e., $\log n$ scale) but differs in the $\log \log n$ scale for the \v{C}ech complex with $k = 1$ in the up-connected case.
\end{abstract}
\date{\today}
\maketitle
\noindent\textit {Key words and phrases.} Random geometric complexes, random hypergraphs, Connectivity, Maximal faces,\\ Phase transition, Poisson convergence.


	\section{Introduction}
	
	Let $\cP_n =\{ X_1,\ldots,X_{N_n}\}$ be a collection of points in the $d$-dimensional torus $U$ with $\{X_i\}_{i \geq 1}$ being a sequence of i.i.d. uniform random variables in $U$ and $N_n$, an independent Poisson random variable with mean $n$. In other words, $\cP_n$ is a homogeneous Poisson point process with intensity $n$ on $U$. A classical model of random graph $G(\cP_n,r)$ (for $r \in (0,\infty)$) introduced by Gilbert in 1961 (\citep{Gilbert61}) called the random geometric graph or Gilbert graph is as follows : The vertex set is $\cP_n$ and $X_i,X_j$ share an edge if $|X_i - X_j| \leq 2r$. Though Gilbert introduced it on the plane, we shall study it on the torus $U$ to avoid boundary effects. This is a common simplification especially when studying sharp thresholds for connectivity properties. A seminal result in the subject was the determining of exact connectivity threshold (\citep{Appel02,Penrose97,Penrose03}). The precise statement of the sharp phase transition result (\citep[Theorem 13.10]{Penrose03}) is that for any sequence $w(n) \to \infty$, the following holds :
	\begin{equation}
	\label{eqn:b0_sharp_threshold}
	\pr{G(\cP_n,r_n) \, \, \mbox{is connected}} \to \begin{cases} 0  \, \, \, & \, \, \, \mbox{if}  \, \, \, n\theta_d2^dr_n^d = \log n  - w(n) \\
	1  \, \,  &  \, \, \, \mbox{if} \, \, \, n\theta_d2^dr_n^d = \log n  + w(n), 
	\end{cases}			   	
	\end{equation} 
	where $\theta_d$ denotes the volume of the unit ball in $\mR^d$. An important step towards the proof of the above result was a similar phase transition result for vanishing of $J_{n,0}$, the number of isolated nodes in  $G(\cP_n,r_n)$ i.e., thresholds for $\pr{J_{n,0} \geq 1}$. That the threshold for vanishing of isolated nodes and threshold for connectivity coincide for a random geometric graph was inspired by a similar phenomenon observed in the case of the Erd\"{o}s-R\'{e}nyi random graphs (\citep{Erdos59}) though the proof in the former case is a lot more involved. 
	
	$J_{n,0}$ is nothing but a $1$-clique and it is natural to wonder if there is a sharp phase transition for vanishing of higher-order cliques and if so, is it related to any higher-dimensional topological phase transitions in random geometric graphs? We denote by $J_{n,k}$ the number of `isolated' $(k+1)$-cliques in $G(\cP_n,r)$, i.e., the number of $(k+1)$-cliques that do not belong to a $(k+2)$-clique. In other words, $J_{n,k}$ is the number of maximal cliques of order $(k+1)$. Due to non-monotonicity of $J_{n,k}$ in $r$, a threshold need not even exist. 
	
	The question of weak/sharp thresholds for higher-order connectivity entails two steps - (1) Determining the weak/sharp threshold for vanishing of `isolated' clique counts and (2) Show that this approximates the weak/sharp threshold of the corresponding notion of `connectivity'. In this article, we shall focus on the first step. One of our results will give thresholds (i.e., $r_n$) for vanishing of $J_{n,k}$ on $G(\cP_n,r_n)$  and for the case of $1 \leq k \leq d$, we shall show that our thresholds are sharp by showing that $J_{n,k}$ converges to a suitable Poisson random variable.  Of course, such thresholds will have implications for the second step too i.e., `connectivity' thresholds. In the next subsection (Section \ref{sec:up_down_connectivity}), we shall discuss in detail about the background literature on these combinatorial notions of connectivity and then preview our results in Section \ref{sec:preview} as well as explaining further connections to existing results. The combinatorial topology notions that we shall need for these two subsections as well as the article are defined rigorously in Section \ref{sec:comb_top_notions}. Our main results are stated in Sections \ref{sec:isol_complex}-\ref{sec:poisson_convg} and we end with proofs in Section \ref{sec:proofs}. 
\vspace*{-0.3cm}
\subsection{Up, Down Connectivity and related literature:}
\label{sec:up_down_connectivity}
	
	A natural higher-dimensional generalization of graphs are simplicial complexes, which in concise terms can be defined as hypergraphs closed under the operation of taking subsets of edges. The precise definition will be given in Section~\ref{sec:comb_top_notions}. The analogous notion of clique counts in simplicial complexes is `face counts'. We provide weak thresholds for vanishing of `isolated' face counts in two models of random geometric complexes - Vietoris-Rips complexes and \v{C}ech complexes. In each of these models, we shall consider two notions of connectivity  - up and down - and hence two notions of `isolation'. We shall define them shortly. 
	
	First, we would like to mention that our question or answer is not without precedence. Specifically, it was shown by Kahle in \citep{Kahle09,Kahle14sharp} that the threshold for vanishing of higher Betti numbers (a notion of higher-order connectivity) was linked to the threshold for vanishing of `isolated' clique counts of {\em Erd\"{o}s-R\'{e}nyi flag/clique complexes}. The $0$th Betti number is nothing but the number of connected components and hence this sharp phase transition result is a generalization of the Erd\"{o}s-R\'{e}nyi result. An earlier generalization of Erd\"{o}s-R\'{e}nyi result for a different model of random complexes called the {\em random $d$-complex} was shown by Linial and Meshulam and later by Meshulam and Wallach \citep{Linial06,Meshulam09}.  Both of these are models that generalize Erd\"{o}s-R\'{e}nyi graphs. We shall not discuss much further about these models of random complexes apart from referring the reader to \citep{Kahle14survey} for more details. The search for a geometric counterpart to the above results is still on despite a significant recent contribution by Bobrowski and Weinberger \citep{Bobrowski15} which we shall discuss later. This is the broader aim towards which we take a step in this article.  Betti numbers represent an algebraic notion of higher-dimensional connectivity and there are other more combinatorial notions of connectivity as we have indicated above and shall discuss now.       
	
	Let $\cK$ be a finite simplicial complex (to be abbreviated as complex in future) and $S_k(\cK)$ denote the set of $k$-faces. The simple notion of connectivity in the graph case generalises to multiple notions of connectivity on  complexes. We shall examine two such notions on two random geometric complexes. Given a complex $\cK$, define the graph of `up-connectivity', $G_k^{\cK,\cU}$ as follows : The vertex set is $S_k(\cK)$ and $\sigma,\tau \in S_k(\cK)$ have an edge if $\sigma \cup \tau \in S_{k+1}(\cK)$. On $\cK$, one can also define the graph of `down-connectivity', $G_k^{\cK,\cD}$ as follows : The vertex set is $S_k(\cK)$ and $\sigma,\tau \in S_k(\cK)$ have an edge if $\sigma \cap \tau \in S_{k-1}(\cK)$. In each of the next four paragraphs, we shall explain four different and unrelated contexts in which `up-connectivity' and `down-connectivity' have been considered. Thus, we hope to convince the reader that these notions of connectivity are worthy of further research not only for their intrinsic challenge and interest but also for their applications.
	
	Having defined a graph, one can naturally consider connected components of the graph, random walk on the graph and the corresponding Laplacian. We shall denote the number of connected components of $G_k^{\cK,\cU}$ and $G_k^{\cK,\cD}$ as $P_k$ ($P$-vector) and $Q_k$ ($Q$-vector) respectively, $k = 0,1,\ldots$.  Since $S_{-1}(\cK) = \emptyset$ by convention, trivially $Q_0 = 0$.  The choice of notation $Q$ has its origins in Q-analysis pioneered by R. Atkins in \citep{Atkin74,Atkin76} to model connectivity of social networks. This has been later developed into a general theory of connectivity of complexes known as combinatorial homotopy or A-homotopy theory. The $Q$-vector plays the role of invariants in this combinatorial homotopy theory and note that $Q_k = 1$ iff $G_k^{\cK,\cD}$ is connected. For more on this combinatorial homotopy theory and its applications, please refer to \citep{Barcelo01,Barcelo05,Kraetzl01}. This is the first context in which the notion of `down-connectivity' is relevant.
	
	Now to the second context. One procedure to construct a complex from a graph $G$ is to define $S_k(\cK)$ to be the set of all $(k+1)$-cliques in $G$. Such a complex is called clique complex of the graph $G$ and we denote it by $\cK(G)$.  Viewing cliques as communities and to investigate overlapping of communities, Derenyi et al. in \citep{Derenyi05} studied percolation on the graph $G_k^{\cK(G),\cD}$ and termed it as clique percolation (though not using the terminology of simplicial complexes). This and further variants of clique percolation on Erd\"{o}s-R\'{e}nyi graphs was studied by \citep{Bollobas09} and the corresponding question for percolation in the up-connectivity graph of random geometric complexes was addressed in \citep{Yogesh13}.  For a survey of this direction of research, see \citep{Palla08}.         
	
	Even though the notion of down-connectivity has implicitly been used in Q-analysis, combinatorial homotopy and clique percolation without stating them explicitly, we shall now reference literature where these terms have appeared explicitly. These are the very recent studies of Laplacians on simplicial complexes (\citep{Horak13,Parzanchevski16,Mukherjee16,Gundert12}), which is the third context where both `up' and `down' connectivity have been studied. Here again, there are two notions of Laplacians - up and down - and as expected they are related to up and down connectivities respectively. Irreducibility of the two Laplacians are related to the connectivity of $G_k^{\cK,\cU}$ and $G_k^{\cK,\cD}$ respectively. As we can observe that there are varied contexts in which the notions of up and down connectivity crop up but barring these few papers on Laplacians of random complexes and face percolation on random complexes, this is very much a fertile terrain. Our results give a lower bound on the thresholds for irreducibility of the two Laplacians and triviality of the $Q$-vector. 
	
	Now the fourth context, which we have touched upon earlier, is the study of Betti numbers of random complexes. This is also the direction of more extensive research on random complexes compared to the  directions in the previous three paragraphs. Betti numbers are an alternate way of quantifying connectivity of complexes. We shall not be formally defining but shall later hint at connection between our results and Betti numbers. We refer the reader to \citep{Edelsbrunner10,Munkres84} for more details on Betti numbers.  In fact, Betti numbers of random complexes has been the main focus of most studies on random complexes. We point the reader to the two surveys \citep{Kahle14survey,Bobrowski14} for details on this growing area lying in the intersection of probability and combinatorial topology. Motivated by applications to topological data analysis (\citep{Carlsson14,Costa12}), random geometric complexes were introduced in \citep{Kahle11} and among other things, upper bounds for thresholds on vanishing of Betti numbers for the Vietoris-Rips and \v{C}ech complexes on Poisson point processes were given. Similar thresholds were later proven for more general stationary point processes in \citep{Yogesh15} and for Poisson point processes on compact, closed manifolds without a boundary in \citep{Bobrowski15}. To briefly allude to applications of threshold results in topological data analysis, we mention that a very weak threshold was used in the pioneering work of \citep{Niyogi08} to find homology of submanifolds from random samples. Since Betti numbers are algebraic quantities and the notions of up-connectivity and down-connectivity are combinatorial, apriori it not obvious why the two need to be related. However, for the random $d$-complex, it was shown that the threshold for up-connectivity was same as the threshold for vanishing of Betti numbers (\citep[Theorem 1.8]{Kahle14}). It is worth repeating that this work is a step towards the geometric counterpart of such a result.  
\vspace*{-0.7cm}
\subsection{Preview : Coarse-scale asymptotics}
\label{sec:preview}
	
	Now, we shall survey some results of relevance to us before stating a few of our results. Our results stated in this section shall be at the coarse-scale as the finer-scale results involve considerably more notation and are postponed to Section \ref{sec:main}. Now onwards, when we refer to random Vietoris-Rips ($\cR(\cP_n,r_n)$) and \v{C}ech complex ($\cC(\cP_n,r_n)$), we refer to these complexes constructed on $\cP_n$, the homogeneous Poisson point process with intensity $n$ on the unit $d$-dimensional torus for some $d \geq 2$ (precisely defined in Definitions \ref{def:rips_complex} and \ref{def:cech_complex}). For the random \v{C}ech complex, a significant contribution refining the afore-mentioned vanishing thresholds appeared recently in \citep{Bobrowski17RSA}. In particular, it was shown that (see \citep[Theorem 5.4]{Bobrowski17RSA})
	for a sequence $w(n) \to \infty$, the following holds : 
	\begin{equation}
	\label{eqn:Omer1}
	\EXP{\beta_k(\cC(\cP_n,r_n))} \to \begin{cases} \infty  \, \, \, & \, \, \, \mbox{if}  \, \, \, n\theta_d r_n^d = \log n + (k-2) \log \log n - w(n) \\
	\beta_k(U)  \, \,  &  \, \, \, \mbox{if} \, \, \, n\theta_d r_n^d = \log n  + k\log \log n + w(n), 
	\end{cases}			   	
	\end{equation} 
	where $\beta_k(U)$ denotes the $k$th Betti number of the $d$-dimensional torus $U$. Apart from the gap between the upper and lower thresholds, this is an extension of \eqref{eqn:b0_sharp_threshold} to higher dimensions. The above threshold result was extended to a phase transition  result for the event $\{\beta_k(\cC(\cP_n,r_n)) \neq \beta_k(U)\}$ at a coarser scale (\citep[Corollary 5.5]{Bobrowski17RSA}) :
	We have for $1 \leq k \leq d-1$ and any $\epsilon \in (0,1)$
	\begin{equation}
	\pr{\beta_k(\cC(\cP_n,r_n)) = \beta_k(U)}\to \begin{cases} 0  \, \, \, & \, \, \, \mbox{if}  \, \, \, n\theta_d r_n^d = (1 - \epsilon)\log n \\
	1  \, \,  &  \, \, \, \mbox{if} \, \, \, n\theta_d r_n^d = (1 + \epsilon)\log n.
	\end{cases}			   	
	\end{equation} 
	The above threshold also corresponds with that of thresholds for complete coverage (\citep{Flatto77,Hall88}) and surprisingly reveals that $\beta_0$ (number of connected components) equals one at a much lower threshold than the other Betti numbers which all vanish ``nearly" together and correspond to the threshold for complete coverage i.e., the event $\{U \subset \cup_{x \in \cP_n}B_x(r_n)\}$. 
	
	As should be obvious by now, the question of connectivity in higher dimensions can be posed in at least three ways (up, down and Betti numbers) and for at least two different models of geometric complexes (Vietoris-Rips and \v{C}ech). Though thresholds for Betti numbers of random \v{C}ech complexes have been partially addressed in \citep{Kahle11,Bobrowski17RSA} but thresholds for other notions of connectivity remain still open. While Betti numbers of \v{C}ech complexes are non-trivial only for $k \leq d-1$ but the question of up and down connectivity are relevant for any $k$.  
	
	The specific geometry of random \v{C}ech complexes enables one to study them via coverage processes as well as Morse theory. Indeed, the key tool in \citep{Bobrowski17RSA} is investigation of the critical points of the  distance function $\rho_n : \mR^d \to \mR,  x \mapsto \min_{1 \leq i \leq N_n} \{ |x-X_i| \}$. Since $\rho_n^{-1}[0,r] = \cup_{X \in \cP_n}B_X(r)$, critical points of $\rho_n$ are related to the Betti numbers of $\cup_{X \in \cP_n}B_X(r)$ via Morse theory and the Betti numbers of the latter are same as that of random \v{C}ech complexes due to nerve theorem. However, both these tools are either unavailable or insufficient for study of other notions of connectivity in the two models of geometric complexes. A possibly more universal approximation for study of connectivity thresholds in random complexes are `isolated' face counts and this is the main reason why we focus on these objects in this article. 
	
	Let $G_{\cdot}^{p,q}$ denote the up and down connected graphs for $q \in \cI_2 := \{\cU,\cD \}$ and
	the \v{C}ech and Vietoris-Rips complexes for $p \in \cI_1 := \{\cC,\cR \}$. In the above notation, $\cC,\cR$ refer to the \Cech and Vietoris-Rips cases respectively and $\cU,\cD$ refer to the up and down connectivity respectively. As a trailer for our results, we state the following coarse scale phase transition result for isolated $k$-faces in the random Vietoris-Rips and \v{C}ech complexes. The constants $m_k^{p,q}$ that appear in the Theorem are defined in Section \ref{sec:comb_top_notions}. 
	\begin{theorem}
		\label{thm:coarse_phase_transition}
		Let $p \in \cI_1, q \in \cI_2$ and $k \geq 1$. Let $J_{n,k}^{p,q}(r_n)$ denote the number of isolated nodes in $G_k^{p,q}(\cP_n,r_n)$. Then the following holds for any $\epsilon \in (0,1)$
		\begin{equation}
		\label{eqn:coarse_phase_transition}
		\pr{J_{n,k}^{p,q}(r_n) = 0}\to \begin{cases} 0  \, \, \, & \, \, \, \mbox{if}  \, \, \, nm_k^{p,q} r_n^d = (1 - \epsilon)\log n \\
		1  \, \,  &  \, \, \, \mbox{if} \, \, \, nm_k^{p,q} r_n^d = (1 + \epsilon)\log n.
		\end{cases}			   	
		\end{equation} 
Indeed, for $nm_k^{p,q} r_n^d = (1 - \epsilon)\log n$, we have that $J_{n,k}^{p,q}(r_n)  \stackrel{P}{\to} \infty$ as $n \to \infty$.
	\end{theorem}

Since $J_n^{p,q}(r_n)$ is a non-monotonic functional in $r$, even presence of a phase transition is not obvious. However, the above theorem shows phase-transition for existence of isolated nodes at a fixed radius $r_n$ in the up/down-connected graphs. However we may also consider existence of isolated nodes for some radii $s \geq r_n$. This is a monotonic event and at the coarse scale has the same threshold for vanishing as existence of isolated nodes. 
	 \begin{theorem}
	 	\label{thm:coarse_phase_transition1}
	 	Let $p \in \cI_1, q \in \cI_2$ and $k \geq 1$. Let $J_{n,k}^{p,q}(r_n)$ denote the number of isolated nodes in $G_k^{p,q}(\cP_n,r_n)$. Then the following holds for any $\epsilon \in (0,1)$
	 	\begin{equation}
	 	\label{eqn:coarse_phase_transition1}
	 	\pr{\cap_{r \geq r_n}\{J_{n,k}^{p,q}(r) = 0\}}\to \begin{cases} 0  \, \, \, & \, \, \, \mbox{if}  \, \, \, nm_k^{p,q} r_n^d = (1 - \epsilon)\log n \\
	 	1  \, \,  &  \, \, \, \mbox{if} \, \, \, nm_k^{p,q} r_n^d = (1 + \epsilon)\log n.
	 	\end{cases}			   	
	 	\end{equation} 
	 \end{theorem}
	We shall shortly see evidence that at a finer scale the  thresholds for the events in Theorems \ref{thm:coarse_phase_transition} and \ref{thm:coarse_phase_transition1} need not coincide at least for the \v{C}ech complex. Considering Theorem \ref{thm:coarse_phase_transition} as the first step towards determining thresholds for up/down-connectivity in both the random geometric complexes, here is the second step. The following theorem shows that whenever `isolated' nodes vanish in the up/down-connected graphs, components of finite but fixed order also vanish. 
	\begin{theorem}
		\label{thm:finite_components}
		Let $p \in \cI_1, q \in \cI_2, k \geq 1$ and $L \geq 1$. Let $J_n^{p,q}(r,L)$ denote the number of components in $G_k^{p,q}(\cP_n,r)$ with exactly $L$ vertices. Then for any $\epsilon > 0$ and $nm_k^{p,q}r_n^d = (1 + \epsilon) \log n$, we have that
		\[ \EXP{J_{n,k}^{p,q}(r_n,L)} \to 0	.\]
	\end{theorem}
Vanishing of isolated nodes in $G^{P,q}_{n,k}(r_n)$ is a necessary condition for connectivity of $G^{P,q}_{n,k}(r_n)$ but as mentioned before, this is also a sufficient condition in many random graph models. It is not completely obvious at this point if this is true even in the models we have considered above. We conjecture it to be true at least in the coarse scale with Theorem 1.3 as a partial eveidence.

More finer phase transition results for expectations and a distributional result inside the critical window for $J_{n,k}^{\cR,\cU}$ with $k \leq d$ are stated in Section \ref{sec:main}. An analogous distributional result is currently unavailable for other models or similar statistics (i.e., Morse critical points) in random geometric complexes. An important tool in obtaining the distributional result for the Vietoris-Rips complex is the purely deterministic geometric Lemma \ref{lem:geometric} and such a result is not available for the \v{C}ech complex or Morse critical points. Since stating these results shall involve considerably more notation and hence are postponed to Section \ref{sec:main}. However, in the special case of $k = 1$, we can state these results with no additional notation as we shall do so now. Note that $J_{n,1}^{\cR,\cD} = J_{n,1}^{\cC,\cD}$. 
	\begin{proposition}
	\label{prop:exp_Jn1}
     Let $J_{n,1}^{p,q}(r_n)$ denote the number of isolated nodes in $G_1^{p,q}(\cP_n,r_n)$.  Let $w_n$ be any real sequence converging to $\infty$ and $nr_n^d \to \infty$ as $n \to \infty$. 
For $(p,q) \in \cI_1 \times \cI_2 \setminus \{(\cC, \cU) \}$ we have 
     \begin{equation}
     \EXP{J_{n,1}^{p,q}(r_n)}\to \begin{cases} \infty  \, \, \, & \, \, \, \mbox{if}  \, \, \, nm_1^{p,q} r_n^d = (\log n - w_n)\\
     0  \, \,  &  \, \, \, \mbox{if} \, \, \, nm_1^{p,q} r_n^d = \log n + w_n,
     \end{cases}			   	
     \end{equation} 
whereas
     \begin{equation}
     \EXP{J_{n,1}^{\cC,\cU}(r_n)}\to \begin{cases} \infty  \, \, \, & \, \, \, \mbox{if}  \, \, \, nm_1^{\cC,\cU} r_n^d = (\log n - \log \log n - w_n)  \\
     0  \, \,  &  \, \, \, \mbox{if} \, \, \, nm_1^{\cC,\cU} r_n^d = \log n - \log \log n + w_n.
     \end{cases}			   	
     \end{equation} 
	\end{proposition}
	%
	
	Notice the difference in the thresholds at the level of expectation between \v{C}ech and Vietoris-Rips complexes. We shall also see that there is a difference between the two complexes when we look at a finer phase transition result corresponding to the event in Theorem \ref{thm:coarse_phase_transition1}. To be more precise, consider the event $\cup_{s \geq r} \{J_{n,1}^{p,q}(s) \geq 1\}$. We denote by $J_{n,k}^{p,q,*}(r)$ the number of $k$-faces (i.e., $k+1$ points) that contribute to the event or in other words, the number of $k$-faces that are `isolated' for some $s \geq r$. This is defined more precisely in \eqref{eqn:J_n^*}. Clearly, all `isolated' faces at radius $r$ are included i.e., $J_{n,k}^{p,q}(r) \leq J_{n,k}^{p,q,*}(r)$ and the latter is non-increasing in $r$. Thus, we have that
$$ \{ J_{n,k}^{p,q,*}(r) = 0 \} = \cap_{s \geq r} \{J_{n,k}^{p,q}(s) = 0\}$$
and we have shown in Theorems \ref{thm:coarse_phase_transition} and \ref{thm:coarse_phase_transition1}	that the thresholds for vanishing of $J_{n,k}^{p,q}(.)$ and $J_{n,k}^{p,q,*}(.)$ coincide at the coarse scale. However, we shall see now that at a finer scale whether they coincide or not depends on the geometry of the complex.
\begin{proposition}
	\label{prop:exp_Jn1s}
	Let $p \in \cI_1, q \in \cI_2$ and $k = 1$. Let $J_{n,1}^{p,q,*}(r_n)$ denote the number of $1$-faces that are isolated in $G_1^{p,q}(\cP_n,s)$ for some $s \geq r_n$.  Let $w_n$ be any real sequence converging to $\infty$ and $nr_n^d \to \infty$ as $n \to \infty$. Then
	\begin{equation}
	\EXP{J_{n,1}^{p,q,*}(r_n)}\to \begin{cases} \infty  \, \, \, & \, \, \, \mbox{if}  \, \, \, nm_1^{p,q} r_n^d = (\log n - w_n) \\
	0  \, \,  &  \, \, \, \mbox{if} \, \, \, nm_1^{p,q} r_n^d = \log n + w_n.
	\end{cases}			   	
	\end{equation} 
\end{proposition}

Note the missing $\log \log n$ factor for the \v{C}ech case i.e., as claimed earlier the thresholds for vanishing of $J_{n,1}^{\cC,\cU}(.)$ and $J_{n,1}^{\cC,\cU,*}(.)$ do not coincide at least at the level of expectations. More importantly, this means that if we choose $w_n \to \infty$ such that  $w_n= o(\log \log n)$ and $nm_1^{\cC,\cU}r_n^d = \log n - \log \log n + w_n$, then
$$ \EXP{J_{n,1}^{\cC,\cU}(r_n)} \to 0, \, \, \, \mbox{but} \, \, \, \EXP{J_{n,1}^{\cC,\cU,*}(r_n)} \to \infty.$$
In other words, for any fixed ``$r \in (\log n - \log \log n, \log n)$" and for $n$ large, we are unlikely to observe an `isolated' face at $r$ but we expect to see large number of them in the interval. 

We shall now explain the consequences of our results for topological phase transitions as well as contrast them with the related results of \citep{Bobrowski17RSA}. While our descriptions are mostly in terms of thresholds for vanishing of isolated nodes but as mentioned before, our tendentious view is that these thresholds are indicators of similar behaviour by the corresponding connectivity thresholds.
\begin{remark}
\label{Remarks1}
\begin{itemize}
\label{list:consequences1}
			
			\item Since $m_k^{\cC,\cU} = \theta_d$ for all $k \geq 1$, we note that at the coarse scale the threshold for vanishing of isolated nodes in $G_k^{\cC,\cU}(\cP_n,r_n)$ matches with that of $\beta_k(\cP_n,r_n)$. And like we pointed out for Betti numbers, $J_{n,0}^{\cC,\cU}$ vanishes much earlier compared to $J_{n,k}^{\cC,\cU}$ for $k \geq 1$ which all vanish ``nearly together".
			
			\item  At a finer scale (see Proposition \ref{prop:convergence_expectation}), the threshold for vanishing of isolated nodes in $G_{.}^{\cC,\cU}$ is analogous to that in \eqref{eqn:Omer1} and in the case of $k =1$ matches exactly with the lower threshold in \eqref{eqn:Omer1} (see Theorem \ref{thm:scaling_r_n_k=1} and Corollary \ref{cor:rate_of_convergence}). In \citep{Kahle14sharp,Linial06,Meshulam09}, it was shown that thresholds for vanishing of isolated nodes in $G_{.}^{.,\cU}$ was same as that of vanishing of homology groups in Erd\"{o}s-R\`{e}nyi-like random complexes. Our results together with \eqref{eqn:Omer1} offer evidence of such a phenomena holding true even for random \v{C}ech complexes.  
			
			\item Perhaps, a little ambitiously one can conjecture that the threshold for vanishing of isolated nodes in $G_k^{\cC,\cU}$ should be $n\theta_dr_n^d = \log n + (k-2) \log \log n$, which corresponds to the lower threshold in \eqref{eqn:Omer1}. Further, in \citep[Corollary 6.3]{Bobrowski17RSA}, it is shown that thresholds for vanishing of critical points of index $k$ are $n\theta_dr_n^d = \log n + (k-1) \log \log n$. Since both index $k$ and $k+1$ critical points are related to $\beta_k(\cC(\cP_n,r_n))$, it is now a moot point as to which of the three thresholds $n\theta_dr_n^d = \log n + l \log \log n$, $l \in \{k-2,k-1,k\}$ is the actual threshold for vanishing of $\beta_k(\cC(\cP_n,r_n))$. 
			
			\item  The four thresholds corresponding to vanishing of isolated nodes in $G_{.}^{\cC,\cU}, G_{.}^{\cC,\cD}, G_{.}^{\cR,\cU}, G_{.}^{\cR,\cD}$ are all different even in the coarser scale, since the corresponding $m_k^{.,.}$'s are different. 
			
			\item Again from the definitions of $m_{.}^{.,\cU}$'s, we can see that vanishing of isolated nodes in $G_{.}^{.,\cU}$ occur for a larger radius than in $G_{.}^{.,\cD}$ even in the coarse scale. Similarly, even in the coarse scale, threshold for vanishing of isolated nodes in $G_{.}^{\cC,.}$ is larger than in $G_{.}^{\cR,.}$.
			
\item Since $m_k^{\cR,.}$ is strictly decreasing in $k$ for $k \leq d$, vanishing of isolated nodes of $G_{k}^{\cR,.}$ happen later than $G_{j}^{\cR,.}$ for $j < k \leq d$ even at the coarse scale. This is in contrast to the scenario for $G_{.}^{\cC,\cU}$ for which $m_{k}^{\cC,\cU} = \theta_d\,$ for all $k \geq 1$.

\end{itemize}
	\end{remark}
	
	We shall now say a few words on our proof methods. While the main tools for obtaining asymptotics for `isolated' faces are the classical Palm theory and Campbell-Mecke formula, the specific geometric analysis pertaining to `isolated faces' differ from that of Morse critical points. Apart from this, we shall see that the non-monotonicity of `isolated' nodes complicates the analysis. In such high-density regimes, the asymptotics are usually determined by the `minimal configuration' contributing to the functional. However, it is very important to understand the behaviour in the neighbourhood of the `minimal configuration'. This we shall see is far from clear for `isolated' faces in contrast to Morse critical points. In special cases such as $k = 1$, we are able to understand the neighbourhood of the `minimal configuration'  to be able to give more detailed results though the techniques vary considerably from case to case. As is evident in our results, such problems are also a matter of scale. Often in coarse scale, we are able to overcome these issues with respect to `the minimal configuration' more easily. 
	
	The connection between Morse critical points and Betti numbers is a deterministic fact whereas the connection between `isolated' faces and Betti numbers arises mainly in random contexts. This is one reason why translating our results to thresholds for Betti numbers is incomplete at the moment. However, for the two discrete models of random simplicial complexes - Erd\"{o}s-R\'enyi clique complexes (\citep{Kahle14sharp}) and random $d$-complexes (\citep{Linial06,Meshulam09}) - using different methods it has been shown that the threshold for vanishing of `isolated' faces in the up-connected graph corresponds to the threshold for vanishing of homology. 
	
In the graph case, closely related to connectivity threshold is the largest edge-weight on a minimal spanning tree. It was shown in \citep{Penrose97} that asymptotically the longest edge on a minimal spanning tree on the complete graph on $\cP_n$ with weights as the Euclidean distance is same as the largest nearest neighbour distance i.e., the smallest radius at which all the vertices are non-isolated in the random geometric graph on $\cP_n$. Simplicial analogues of spanning tree are called as {\em spanning acycles} and their behaviour on randomly weighted complexes have been investigated in \citep{hiraoka2015minimum,Skraba2017}. In particular, see \citep[Section 3]{Skraba2017} for relations between spanning acycles, Betti numbers and `isolated' faces. We see our work as another step towards establishing such relations for random geometric complexes.	
	
To end the introduction, we shall point a few more directions in which our work could be extended apart from the natural program of computing exact thresholds and investigating the critical window. One important problem would be to investigate higher-dimensional or distributional analogues of  Propositions \ref{prop:exp_Jn1} and \ref{prop:exp_Jn1s}. One could also consider geometric complexes on compact Riemannian manifolds and Poisson point processes with non-uniform densities. For ideas on the former extension, see \citep{Bobrowski17RSA,Bobrowski17} and for the latter see \citep{Gupta2010,Iyer2012,Penrose03,Owada15,Owada16}.  

\section{Main Results}
\label{sec:main}

\subsection{Some combinatorial topology notions :}
\label{sec:comb_top_notions}

A subset $\cK \subset 2^{\cX}$ for a finite point-set $\cX$ is said to be an \emph{abstract simplicial complex} (abbreviated as complex in future) if $A \in \cK$ and $B \subset A$ implies that $B \in \cK$. The elements of $\cK$ are called faces or simplices and the dimension $dim \, \sigma$ of a face $\sigma$ is $|\sigma| -1$. We shall denote a $k$-face by $[v_1,\ldots,v_{k+1}]$.  The maximal faces (faces that are not included in any other faces) are called {\em facets}. By convention, $\emptyset \in \cK$ and $dim(\emptyset) = -1$. The collection of $k$-faces is denoted by $S_k(\cK)$ and the $k$-skeleton of $\cK$ is the complex $\cK^k := \cup_{i=-1}^k S_i(K)$. A complex is said to be {\em pure} if all maximal faces have the same dimension. Note that $\cK^1$ is nothing but a graph. 

Denote by $B_x(r)$ a closed ball of radius $r$ centered at $x$. $| \cdot |$ will denote the cardinality of a finite set as well as the Lebesgue measure and $\| \cdot \|$ is the Euclidean norm on $\mR^d$. There are two types of complexes defined on point processes. These complexes are the Vietoris-Rips complex and the \Cech complex which we now define. Let $U = [0,1]^d$ be equipped with the {\em toroidal metric} i.e.,
$$ d(x,y) = \inf \{\|x-y+z\| : z \in \mZ^d\}, \, \, x,y \in U.$$
For $\x = (x_1,\ldots,x_{k+1}) \in \mR^{d(k+1)}$, let $B_{\x}(r) = \bigcup_{i=1}^{k+1} B_{x_i}(r)$, $h(\x) = h(x_1,\ldots,x_{k+1})$ for $h \: \mR^{d(k+1)} \to \mR$ and $\md \x = \md x_1 \ldots \md x_{k+1}$. Let $\mathbf{1} = (1,\ldots,1)$. Let $\cX$ be a finite set in $U$ and we shall use $\cX^{(k)}$ to denote the set of $k$-tuples of distinct points in $\cX$.

\begin{definition}[Vietoris-Rips complex]
	\label{def:rips_complex}
	The abstract simplicial complex $\cR(\cX,r)$ constructed as below is called the \VR complex associated to $\cX$ and $r$.
	\begin{enumerate}
		\item The $0$-simplices of  $\cR(\cX,r)$    are the points in $\cX$.
		\item A $k$-simplex, or $k$-dimensional `face',  $\sigma=[x_{i_1},\ldots,x_{i_{k+1}}]$ is in $\cR(\cX,r)$ if $B_{x_{i_j}}(r) \cap B_{x_{i_m}}(r) \ne \emptyset$ for every $1 \le j < m \le k+1$ and where $(x_{i_1},\ldots,x_{i_{k+1}}) \in \cX^{(k+1)}$.
	\end{enumerate}
\end{definition}
\begin{definition}[\Cech complex]
	\label{def:cech_complex}
	The abstract simplicial complex  $\cC(\cX,r)$  constructed as below is called the \Cech complex associated to $\cX$ and $r$.
	\begin{enumerate}
		\item The $0$-simplices  of  $\cC(\cX,r)$   are the points in $\cX$,
		\item A $k$-simplex, or $k$-dimensional `face',  $\sigma=[x_{i_1},\ldots,x_{i_{k+1}}]$ is in $\cC(\cX,r)$ if $\bigcap_{j=1}^{k+1} B_{x_{i_j}}(r) \ne  \emptyset$ and where $(x_{i_1},\ldots,x_{i_{k+1}}) \in \cX^{(k+1)}$.
	\end{enumerate} 	
\end{definition}
Observe that the faces of a Vietoris-Rips complex are nothing but cliques of a random geometric graph and the $1$-skeletons of both the Vietoris-Rips and \Cech complexes coincide with the random geometric graph.

\begin{figure}[!htbp]
\centering
\includegraphics[width=3in,height=2.5in]{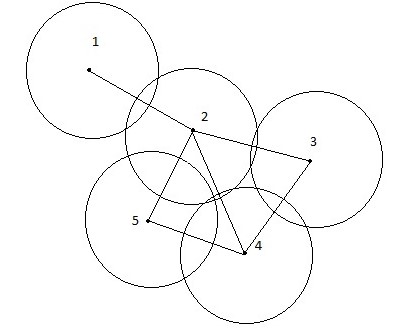}
\caption{$\{[1,2],[2,3],[3,4],[2,4,5]\}$ are the maximal faces of the \Cech complex on the point set $\cX = \{1.\ldots,5\}$ whereas $\{[1,2],[2,3,4],[2,4,5]\}$ are the maximal faces of the Vietoris-Rips complex on $\cX$.}
\label{fig_complex}
\end{figure}

The functionals that we study in this paper are the isolated simplex counts in the \Cech and the Vietoris-Rips complexes. In addition there are two notions of connectivity, up and down which in turn determines what constitutes an isolated simplex. For any $k \geq 0$, let $S_{k}(\cX,r)$ be the collection of all $k$-simplices of the \Cech complex on $\cX$. Consider the graph $\Gcu(\cX,r)$ with vertex set $S_{k}(\cX,r)$ and with an edge between any two elements $\sg_1, \sg_2 \in S_{k}(\cX,r)$ provided they are up-connected, that is, $\sg_1 \cup \sg_2 \in S_{k+1}(\cX,r)$. Similarly for any $k \geq 1$, we can define the graphs $\Gcd(\cX,r)$ with edges between elements $\sg_1, \sg_2 \in S_{k}(\cX,r)$ that are down connected, that is, $\sg_1 \cap \sg_2 \in S_{k-1}(\cX,r)$. The graphs $\Gru(\cX,r)$ and $\Grd(\cX,r)$ are defined similarly by taking $S_k(\cX,r)$ to be the collection of all $k$-simplices in the Vietoris-Rips complex. See the below figure~\ref{fig_complex} for an illustration of the two complexes and their maximal faces. 
%

Fix $k \geq 2$. For $\x =(x_1, \ldots ,x_{k+1}) \in \mR^{d(k+1)}$ and $r > 0$ define the functions
\begin{eqnarray}
\hc_k(\x,r) & = & \1[\bigcap_{j=1}^{k+1} B_{x_{j}}(r) \ne  \emptyset], \nn \\
\hr_k(\x,r) & = & \prod_{1 \leq j < m \leq k+1}\1[B_{x_j}(r) \cap B_{x_m}(r) \ne \emptyset]. 
\lab{def:isomorphism}
\end{eqnarray}
Recall that $\cI_1 = \{\cC,\cR\}$ and $\cI_2 = \{\cU,\cD\}.$ We will denote by $h^p_k(\x)$ the function $h^p_k(\x,1)$ for $p \in \cI_1, k \geq 1$. For $x \in \mR^d$ and $r,s > 0$ let $\lb B_x(r) \rb^{(s)} = B_x(r) \oplus B_O(s)$ be the closed $s-$neighbourhood of the ball $B_x(r)$ and $O$ denotes the origin. For $x \in \mR^{d(k+1)}$ and $r,s > 0$ define the set-valued functions
\begin{eqnarray}
\Qcu_k(\x,r,s) & = & \lb \cap_{j=1}^{k+1} B_{x_j}(r) \rb^{(s)} \nn \\
\Qcd_k(\x,r,s) & = & \cup_{i=1}^{k+1} \lb \cap_{j=1, j \ne i}^{k+1} B_{x_j}(r) \rb^{(s)} \nn \\
\Qru_k(\x,r,s) & = & \cap_{j=1}^{k+1} \lb B_{x_j}(r) \rb^{(s)} = \cap_{j=1}^{k+1} B_{x_j}(r+s) \nn \\
\Qrd_k(\x,r,s) & = & \cup_{i=1}^{k+1} \lb \cap_{j=1, j \ne i}^{k+1} \lb B_{x_j}(r) \rb^{(s)} \rb.
\lab{def:Qfns}
\end{eqnarray}
We will also use the abbreviated forms 
\begin{equation}
Q^{p,q}_k(\x,r) = Q_k^{p,q}(\x,r,r) \qquad  \mbox{and} \qquad Q^{p,q}_k(\x) = Q^{p,q}_k(\x,1), \qquad p \in \cI_1, \;\; q \in \cI_2.
\lab{eqn:Qspec}
\end{equation}
Let
\begin{equation}
A^p_k := \{ \y \in (\mR^d)^{k} : h_k^p(O,\y) = 1 \},
\lab{eqn:A_p}
\end{equation}
be the set of configurations that form a $k$-simplex with the origin. Clearly $A^p_k \subset B_O(2)^{k}$ and since our complexes are defined using closed balls, $A^p_k$ is compact. Let 
\begin{equation}
m_k^{p,q} = \inf\{|Q_k^{p,q}(O,\y)|: \y \in A_k^p \}, \qquad \mbox{and} \qquad 
M_k^{p,q} = \sup\{|Q_k^{p,q}(O,\y)|: \y \in A_k^p \}.
\end{equation}
Often, the choice of $k$ is clear from the context and hence we shall suppress it. It is easy to see that $m_k^{\cC,\cU} = \th_d$ and $m_k^{\cC,\cD} = 2^d \th_d$ for all $k \geq 1$. The former occurs for a configuration of points that are as far apart as possible such that the common intersection of balls is a single point and the latter happens when all the points coincide. Note that $m_k^{\cC,q}$ does not depend on $k$. For $k \leq d$, $m_k^{\cR,\cU}$ is the volume of the lens of intersection of $(k+1)$ balls of radius $2$ that can be placed in $\mR^d$ so that their centers are exactly at distance $2$ from each other.

\subsection{Expectation Asymptotics for Isolated face Counts}
\label{sec:isol_complex}
In this section we present a radius regime under which the expected isolated complex count stabilizes in the limit. This regime involves a parameter sequence whose asymptotic behavior is described. Later part of the section contains a result on the rate of convergence of this parameter sequence which has some interesting implications that will be discussed.
\begin{definition}[Isolated face counts]
\label{def:isol_complex_counts}
Let $k \geq 1$. For $p \in \cI_1$ and $q \in \cI_2$ the number of isolated  simplices in the graph $G^{p,q}(\cX,r)$ is defined as
\begin{equation}
J^{p,q}(\cX,r) = J^{p,q}_k(\cX,r) := \frac{1}{(k+1)!}\sum_{\x \in \cX^{(k+1)}}h^p(\x,r) \1[\cX \cap Q^{p,q}(\x,r) \equiv \x].
\lab{def:J}
\end{equation}
\end{definition}

For example, $J^{\cR,\cU}(\cX,r)$ counts the number of maximal $(k+1)$-cliques in the random geometric graph. As a more concrete example, in Figure \ref{fig_complex}, we have that $J_1^{\cR,\cD}(\cX,r) =J_1^{\cC,\cD}(\cX,r) = 0$, $J_1^{\cR,\cU}(\cX,r) = 1, J_1^{\cC,\cU}(\cX,r) = 3.$ Let $\cP_n$ be a Poisson point process with intensity $n\1_{U}(.)$ where $U = [0,1]^d$ is the unit cube with the toroidal metric. We will denote $J^{p,q}(\cP_n,r)$ by $J_n^{p,q}(r)$ and if we wish to emphasize the dependence on $k$, we shall denote by $J_{n,k}^{p,q}(r)$. We shall often drop the subscript $k$ from our other notations such as $A^p_k,M^{p,q}_k, m^{p,q}_k$. Our first result is on the radius regime $r_n$ that stabilizes the expected number of isolated simplices in the connectivity regime.
For any $k \geq 1$, $\alpha \in \mR$, $p \in \cI_1$ and $q \in \cI_2$, $c > 0$ define the sequence of radial functions $\{r_n^{p,q}(c)\}_{n \geq n_0}$ as
\begin{equation} 
r_n^{p,q}(c)  = \lb \f{\log n + k \log \log n + \log |A^p| + \alpha -  k \log m^{p,q} - \log (k+1)!}{n c} \rb^{\f{1}{d}},
\label{def:stabilising_radii}
\end{equation}
where $n_0$ is defined so that for all $n \geq n_0$, $r_n > 0$. Note that $n_0$ does not depend on $c$. 
\begin{proposition}
\label{prop:convergence_expectation}
Let $k \geq 1$, $\alpha \in \mR$, $p \in \cI_1, q \in \cI_2$ and $r_n^{p,q}(c)$ be as defined in (\ref{def:stabilising_radii}). Then there exists a sequence $\{c_n^{p,q}\}_{n \geq 1}  \subset [m^{p,q},M^{p,q}]$ such that $c_n^{p,q} \to m^{p,q}$ and
\begin{equation} 
\EXP{J_n^{p,q}(r_n^{p,q}(c_n^{p,q}))} \to e^{-\alpha}, \qquad \mbox{as } n \to \infty.
\end{equation}
\end{proposition}
Though the constant $c_n^{p,q} \to m^{p,q}$, one cannot replace $c_n^{p,q}$ by $m^{p,q}$ in Proposition~\ref{prop:convergence_expectation} as shown by the following result.
\begin{proposition}
\label{prop:convergence_expectation_to_zero}
Let $k \geq 1$, $\alpha \in \mR$, $p \in \cI_1, q \in \cI_2$ and $r_n$ be such that $nm^{p,q}r_n^d = \log n + k \log \log n + w^1_n$ where $w^1_n$ is a sequence bounded from below i.e.,$\liminf_{n \to \infty} w^1_n > - \infty$. Then 
\begin{equation}
\EXP{J_n^{p,q}(r_n^{p,q}(m^{p,q}))} \to 0, \qquad \mbox{as } n \to \infty.
\end{equation}
\end{proposition}
It would be desirable to obtain precise estimates on the constant $c_n^{p,q}$ in Proposition~\ref{prop:convergence_expectation}. We explore this and a few other results for the case $k = 1$ in the next subsection. 
\subsubsection{The 1-simplex. }

For the case $k=1$ we derive a result on the rate of convergence of the sequence $c_n^{p,q}$ and discuss some interesting implications of these results. The maximal face in a $1$-complex is an edge. An edge is isolated in the sense of up-connectivity provided it is not part of a triangle i.e., a $2$-simplex. For three nodes to form a triangle in the Vietoris-Rips complex, balls centered at these three vertices must have non-trivial pairwise intersection. In the \Cech case the balls centered at these vertices must have a common intersection. Down connectivity for both the Rips and the \Cech cases are identical. An edge is isolated in the down sense if it does not share a vertex with any other edge.

The following theorem shows that the rate of convergence of $c_n^{p,q}$ to $m^{p,q}$ alters the coefficient of the $\log \log n$ term in the numerator of $r_n^{p,q}(c_n^{p,q})$ as in \eqref{def:stabilising_radii} if we express it in terms of $m^{p,q}$ instead of $c_n^{p,q}$.
\begin{theorem}
Let $a^{\cC,\cU} = 2$ and $a^{p,q} = 1$ for $(p,q) \in \cI_1 \times \cI_2 \setminus (\cC,\cU)$. Let $r_n^{p,q}(c_n^{p,q})$, $p \in \cI_1$ and $q \in \cI_2$ be the sequence as in Proposition~\ref{prop:convergence_expectation} for which the expected number of isolated edges, $\EXP{J_{n,1}^{p,q}(r_n^{p,q}(c_n^{p,q}))} \to e^{-\alpha}$ with $k=1$. Then 
\begin{equation}
n \lb r_n^{p,q}(c_n^{p,q}) \rb ^d m^{p,q} - \log n - (1 - a^{p,q}) \log \log n
\lab{eqn:scaling_r_n_k=1}
\end{equation}
is a bounded sequence.
\label{thm:scaling_r_n_k=1}
\end{theorem}
Though the statement of the above theorem covers all the cases via a single equation, the estimates in different cases require somewhat different ideas. By substituting for $r_n^{p,q}(c_n^{p,q})$ from (\ref{def:stabilising_radii}) in (\ref{eqn:scaling_r_n_k=1}) yields the following result on the rate of convergence of $c_n^{p,q}$ to $m^{p,q}$.
\begin{corollary}
\label{cor:rate_of_convergence}
For $p \in \cI_1$, $q \in \cI_2$, $\al \in \mR$ and $k=1$, let $c_n^{p,q}$ be the sequence as in Proposition~\ref{prop:convergence_expectation} for which the expected number of isolated edges, $\EXP{J_{n,1}^{p,q}(r_n^{p,q}(c_n^{p,q}))} \to e^{-\alpha}$ and let $a_{p,q}$ be as defined in Theorem \ref{thm:scaling_r_n_k=1}. Then,
$$ \lo \lb \f{c_n^{p,q}}{m^{p,q}} - 1 \rb \f{\log n}{\log \log n} = a_{p,q}.$$
\end{corollary}
It is somewhat surprising that $a^{\cC,\cU} = 2$ in the contrast to the other cases and the implication of this for the threshold for vanishing of isolated faces has been mentioned in Remark \ref{Remarks1}.
\subsection{Monotonic vanishing of Isolated faces :}
\label{sec:mon_vanishing}
	
Let  $p \in \cI_1$ and $q \in \cI_2$. Given $X_1,\ldots,X_{k+1}$, define \\ $R^{p}(X_1,\ldots,X_{k+1}) := \inf \{r : h^{p}((X_1,\ldots,X_{k+1}),r) = 1\}$. For the Vietoris-Rips complex, we have that $2R^{\cR}(X_1,\ldots,X_{k+1}) = \max_{i \neq j}|X_i-X_j|.$ When the $(k+1)$-tuple $X_1,\ldots,X_{k+1}$ is clear, we shall simply use $R^p$ instead of $R^p(X_1,\ldots,X_{k+1})$. Define {\em the number of isolated faces at $r$ and above} as follows : 
	\begin{equation}
	\label{eqn:J_n^*}
	J^{p,q,*}_k(\cX,r) := J^{p,q}_k(\cX,r) + \sum_{(X_1,\ldots,X_{k+1}) \in \cX^{(k+1)}}\1[R^p > r] \1[\cX \cap Q^{p,q}((X_1,\ldots,X_{k+1}),R^p) = \{X_1,\ldots,X_{k+1}\}],
	\end{equation}
	where $Q^{p,q}(.,.)$'s are defined in \eqref{eqn:Qspec}. By definition, it is clear that $J_k^{p,q,*}(\cX,r) \geq J_k^{p,q}(\cX,r)$. We have already seen coarse-scale thresholds for vanishing of $J_{n,k}^{p,q,*} := J^{p,q}_k(\cP_n,r_n), k \geq 1$ in Theorem \ref{thm:coarse_phase_transition1} and a more finer threshold for vanishing of $J_{n,1}^{p,q,*}$ in Proposition \ref{prop:exp_Jn1s}.  We now give a finer upper bound for the threshold for vanishing of $J_k^{p,q,*}(\cP_n,r)$ for all $k \geq 1$.
	\begin{proposition}
    \label{prop:exp_mon_threshold}
    Fix $k \geq 1, p \in \cI_1, q \in \cI_2$ and let $J_{n,k}^{p,q,*} := J^{p,q}_k(\cP_n,r_n)$ where       $nm_k^{p,q}r_n^d = \log n + k \log \log n + w^1_n$ for some sequence $w^1_n$ bounded from below i.e., $\liminf_{n \to \infty} w_n^1 > - \infty.$ Then we have that
     \[ \EXP{J_{n,k}^{p,q,*}(r_n)} \to 0, \qquad \mbox{as } n \to \infty.\]
	\end{proposition}
	%
	 
\subsection{Poisson convergence for isolated Vietoris-Rips simplices under up-connectivity}
\label{sec:poisson_convg}

Our next result is a weak convergence result for the number of isolated simplices in the Vietoris-Rips complex.
\begin{theorem} Let $\al \in \mR$, $d \geq 2, 0 \leq k \leq d$ and $J_n = J^{\cR,\cU}_n(r_n^{\cR,\cU}(c_n))$ where $\{c_n = c_n^{\cR,\cU}\}_{n \geq 1}$ is the sequence as in Proposition~\ref{prop:convergence_expectation} i.e., $\EXP{J_n} \to e^{-\alpha}$. Then the number of isolated $k$-simplices $J_n$ converges in distribution to a Poisson random variable with mean $e^{-\al}$.
\label{thm:rips-Poisson-cgs}
\end{theorem}
The above distributional result extends to finite connected components in $G^{\cR,\cU}(\cP_n,r_n^{\cR,\cU}(c_n))$. Recall that $J_n^{\cR,\cU}(r,\ell)$ denotes the number of components in $G_k^{\cR,\cU}(\cP_n,r)$ with $\ell$ vertices (see Theorem \ref{thm:finite_components}).
\begin{theorem} Let $\al \in \mR$, $d \geq 2$ and $0 \leq k \leq d$. Suppose that the $\{c_n = c_n^{\cR,\cU}\}_{n \geq 1}$ is the sequence as in Proposition~\ref{prop:convergence_expectation} i.e., $\EXP{J^{\cR,\cU}_n(r_n^{\cR,\cU}(c_n))} \to e^{-\alpha}$. Let $L \geq 1$. Then $\sum_{l=1}^LJ^{\cR,\cU}_n(r_n^{\cR,\cU}(c_n),l)$ converges in distribution to a Poisson random variable with mean $e^{-\al}$.
\label{thm:rips-Poisson-cgs_finite}
\end{theorem}
\section{Proofs}
\label{sec:proofs}
In what follows, $C_1,C_2, \ldots$ will denote finite constants whose values will change from place to place.
\subsection{Proofs of Results in Section \ref{sec:isol_complex}}
\hspace{1cm} \\ \\
{\bf Proof of Proposition~\ref{prop:convergence_expectation}. } Fix $p \in \cI_1$ and $q \in \cI_2$. We will drop the superscripts $p,q$ for the rest of the proof since the proof is identical in all the four cases. For $\x = (x_1,\ldots,x_{k+1})$, let $h_n(\x) = h(\x,r_n)$ and $Q_n(\x) = Q(\x,r_n)$ be as defined in \eqref{def:isomorphism} and \eqref{def:Qfns}. Without loss of generality, we can assume that $n \geq n_0$, where $n_0$ is as chosen in \eqref{def:stabilising_radii}.

For any sequence of radial functions $r_n$, by the Campbell-Mecke formula (see \cite[Theorem 1.6]{Penrose03}) we have
\begin{equation}
\label{eqn:exp_jn1}
\EXP{J_n(r_n)} = \frac{n^{k+1}}{(k+1)!}\int_{U^{k+1}} h_n(\x) e^{-n |Q_n(\x)|} \md \x .
\end{equation}
Setting $x_i = x_1 + r_n y_i$, $i = 2, \ldots ,k+1,$ in (\ref{eqn:exp_jn1}) and using the fact that $|Q_n(\x)| = r_n^d |Q(\x)|$ we get 
\begin{equation}
\label{eqn:exp_jn2}
\EXP{J_n(r_n)} = \frac{n(n r_n^d)^{k}}{(k+1)!}\int_{U \times ((r_n)^{-1}(U - x_1))^k}
h(x_1, \y) e^{- n r_n^d |Q(x_1,\y)|} \md x_1 \md \y ,
\end{equation}
where $\y = (y_2,\ldots,y_{k+1})$. Let $n_1 \geq n_0$ be such that for all $n \geq n_1$ such that $\cup_{x \in U}B_x(1) \subset (r_n)^{-1}(U - x_1))$ for all $x_1 \in U$. Since the metric is toroidal, we obtain for all $n \geq n_1$ that
\begin{equation}
\label{eqn:exp_jn3}
\EXP{J_n(r_n)} = \frac{n(n r_n^d)^{k}}{(k+1)!}\int_{A} e^{- n r_n^d |Q(O,\y)|} \md \y ,
\end{equation}
where $A = A^p$ is as defined in (\ref{eqn:A_p}). 

For any $c \in \mR_+$ and $\alpha \in \mR$, let $r_n(c)$ be as defined in (\ref{def:stabilising_radii}). The function 
\[ c \mapsto  \int_A e^{- n r_n(c)^d |Q(O,\y)|} \md \y \]
is continuous in $c$, and for all $n > n_1$, tends to $0$ as $c \to 0$ and tends to $|A|$ as $c \to \infty$. 
Since $e^{-n r_n(c)^d c} \in (0,1)$ does not depend on $c$, by the intermediate value theorem there exists a sequence $c_n$ such that
\begin{equation}
\label{eqn:c_n}
\int_A e^{- n r_n(c_n)^d |Q(O,\y)|} \md \y = |A| e^{-n r_n(c_n)^d c_n}.
\end{equation}
With the above choice of $c_n$,  we obtain from (\ref{eqn:exp_jn3})-(\ref{eqn:c_n}) that
\begin{eqnarray}
\EXP{J_n(r_n(c_n))}  & = & \frac{n(nr_n(c_n)^d)^{k}}{(k+1)!}|A|e^{-n r_n(c_n)^d c_n} \no \\
		& = & \frac{n(nr_n(c_n)^d)^{k}}{(k+1)!}|A|e^{-(\log n + k \log \log n + \log |A| + \alpha -  k \log m - \log (k+1)!)} \no \\ 
		& = & \left( \frac{n r_n(c_n)^d m}{\log n} \right)^{k} e^{-\alpha} \no \\
		& = & \left( \frac{\log n + k \log \log n + \log |A| + \alpha - k \log m - \log (k+1)!}{\log n} \right)^{k} \left(\frac{m}{c_n} \right)^{k} e^{-\alpha}. 
 \label{eqn:exp_convergence}
	\end{eqnarray}
Thus the proof is complete provided we show that $c_n \in (m,M)$ and $c_n \to m$. From the definition of $m,M$ and the fact that the function $Q(O,\y)$ achieves these values only on a set of zero measure, we have
\[ |A| e^{-n r_n(c)^d M} < \int_A e^{- n r_n(c)^d |Q(O,\y)|} \md \y < |A| e^{-n r_n(c)^d m}, \]
which together with (\ref{eqn:c_n}) implies that $c_n \in (m,M)$.  

Suppose $\ls c_n > m$. Then we can choose $m_1 > m$, $\ep > 0$ and a subsequence $\{n_j\}_{j \geq 1}$ such that $c_{n_j} (1 - \ep) > m_1$ for all $j \geq 1$. From (\ref{eqn:exp_jn3}) we derive by calculations similar to the one in (\ref{eqn:exp_convergence}) that
\begin{eqnarray*}
\lefteqn{\EXP{J_{n_j}(r_{n_j}(c_{n_j}))}  \geq  \frac{n_j (n_j r_{n_j}^d)^{k}}{(k+1)!}\int_{A\cap\{Q(O,\y) < m_1\}} e^{- n_j r_{n_j}^d |Q(O,\y)|} \md \y} \no \\
 & \geq & C n_j (n_j r_{n_j}^d)^{k} \exp\lb - \f{m_1}{c_{n_j}} \lb \log n_j + k \log \log n_j + \log |A| + \alpha -  k \log m - \log (k+1)! \rb \rb \no \\
 & \geq &  C n_j (n_j r_{n_j}^d)^{k} \exp\lb - (1 - \ep) \lb \log n_j + k \log \log n_j + \log |A| + \alpha -  k \log m - \log (k+1)! \rb \rb \to \infty.
\end{eqnarray*}
This contradicts the fact that from \eqref{eqn:exp_convergence} and $c_n > m$, we must have $\limsup_{n \to \infty}\EXP{J_n(r_n(c_n))} \leq e^{-\alpha}$. Hence $\ls c_n = m$ and since $c_n \in (m,M)$, it follows that $\lo c_n = m$. \qed

{\bf Proof of Proposition~\ref{prop:convergence_expectation_to_zero}. } Dropping the superscripts $p,q$ we have from (\ref{eqn:exp_jn3}) that
\begin{eqnarray}
\EXP{J_n(r_n(m))} & = & \frac{n(nr_n^d)^{k}}{(k+1)!} e^{-nmr_n^d} \int_{A^p} e^{- n r_n^d (|Q(O,\y)| - m)} \md \y \no \\
 & \leq & C_1 e^{-w_n^1}\int_{A^p} e^{- n r_n^d (|Q(O,\y)| - m)} \md \y \to 0, \, \, \, n \to \infty.
\end{eqnarray}
The convergence above follows by the bounded convergence theorem because $A^p$ is compact, $|Q(O,\y)| \geq m$ on $A^p$ and $e^{- n r_n^d (|Q(O,\y)| - m)} \to 0$ almost surely on $A^p$. \qed 

{\bf Proof of Theorem~\ref{thm:scaling_r_n_k=1}. } Since $k=1$ we need to consider only three cases as isolated simplex counts for `down-connectivity' in both the Vietoris-Rips and \v{C}ech complexes are identical for $k=1$.

{\it Case 1. } We first consider the case $p = \cR, q = \cU$. Since $k=1$, the function $Q^{\cR,\cU}(O,y)$ in (\ref{eqn:exp_jn3}) equals$|B_{O}(2) \cap B{y}(2)|$, where $y \in A = B_{O}(2)$. Substituting in (\ref{eqn:exp_jn3}) and changing to polar coordinates we get
\begin{equation}
\EXP{J_n(r)} = n (nr^d) \th_d \int_0^2 s^{d-1} e^{-n r^d |B_O(2) \cap B_{se_1}(2)|} \md s, 
\lab{eqn:k=1_expectation_rips_up}
\end{equation}
where $e_1 = (1,0, \ldots ,0) \in \mR^d$ is the unit vector along the first coordinate axis. The volume of the lens of intersection $B_O(2) \cap B_{s e_1}(2)$ equals $2^d \th_d \eta(s)$ (see \citep[(7.5)]{Goldstein10} and \citep[(6)]{Moran73}) where 
\begin{equation}
\eta(s) =  1 - \f{\th_{d-1}}{\th_d} \int_0^{s/2} \lb 1 - \f{t^2}{4} \rb ^{\f{d-1}{d}}, \qquad 0 \leq s \leq 2.
\lab{eqn:vol_lens}
\end{equation}
Substituting from (\ref{eqn:vol_lens}) in (\ref{eqn:k=1_expectation_rips_up}) we get
\begin{equation}
\EXP{J_n(r)} = n (nr^d) \th_d e^{-nr^d \th_d 2^d \eta(2)} \int_0^2 s^{d-1} e^{-n r^d \th_d 2^d (\eta(s) - \eta(2))} \md s.
\lab{eqn:k=1_expectation_rips_up2}
\end{equation}
For $0 \leq t \leq 1$ we have $\f{3}{4} \leq \lb 1 - \f{t^2}{4} \rb \leq 1$ and hence
\begin{equation}
\f{\th_{d-1}}{\th_d} \lb \f{3}{4} \rb ^{\f{d-1}{d}} \lb 1 - \f{s}{2} \rb \leq \eta(s) - \eta(2) \leq \f{\th_{d-1}}{\th_d}  \lb 1 - \f{s}{2} \rb.
\lab{eqn:bound_lens}
\end{equation}
Using the lower bound for $(\eta(s) - \eta(2))$ from (\ref{eqn:bound_lens}) in (\ref{eqn:k=1_expectation_rips_up2}) and noting that $m = m^{\cR,\cU} = \th_d 2^d \eta(2)$, we obtain
\begin{equation}
\EXP{J_n(r)} \leq n (nr^d) \th_d e^{-nr^d m} \int_0^2 s^{d-1} e^{-n r^d 2^d \th_{d-1} \lb \f{3}{4} \rb ^{\f{d-1}{d}} \lb 1 - \f{s}{2} \rb } \md s.
\end{equation}
Let $a = 2^d \lb \f{3}{4} \rb ^{\f{d-1}{d}} \th_{d-1}.$ Making the change of variable $u = n r^d 2^d \th_{d-1} \lb \f{3}{4} \rb ^{\f{d-1}{d}} \lb 1 - \f{s}{2} \rb$ and replacing $r$ by $r_n$ we get
\begin{equation*}
\EXP{J_n(r_n)} \leq C_2 n  e^{-n r_n^d m} \int_0^{a n r_n^d} \lb 1 - \f{u}{a n r_n^d} \rb ^{d-1} e^{-u} du.
\end{equation*}
If $n r_n^d \to \infty$ as $n \to \infty$, then 
\[ \int_0^{a n r_n^d} \lb 1 - \f{u}{a n r_n^d} \rb ^{d-1} e^{-u} du \to 1, \]
and hence
\begin{equation}
\EXP{J_n(r_n)} \leq e^{- n r_n^d m + \log n + C_3}.
\lab{eqn:k=1_expectation_rips_up4}
\end{equation}
Since $\EXP{J_n(r^{\cR,\cU}_n(c_n^{\cR,\cU}))} \to e^{-\al} \in (0,\infty)$, we must have
\[ n (r^{\cR,\cU}_n(c_n^{\cR,\cU}))^d m_2 - \log n \leq C_4. \]
Similarly using the upper bound for $(\eta(s) - \eta(2))$ from (\ref{eqn:bound_lens}) in (\ref{eqn:k=1_expectation_rips_up2}) and proceeding as above will yield
\begin{equation}
\EXP{J_n(r_n)} \geq e^{- n r_n^d m + \log n + C_5},
\label{eqn:k=1_expectation_rips_up5}
\end{equation}
and again using the fact that $\EXP{J_n(r^{\cR,\cU}_n(c_n^{\cR,\cU}))} \to e^{-\al}$, we obtain
\[ n (r^{\cR,\cU}_n(c_n^{\cR,\cU}))^d m_2 - \log n \geq C_6. \]
{\it Case 2. } Let $p=\cC, q = \cU$. For $k=1$, $Q^{\cC,\cU}(O,y) = V_d \lb \lb B_O(1) \cap B_y(1) \rb \oplus B_O(1) \rb$ where $V_d$ denotes the volume in $\mR^d$. Substituting in (\ref{eqn:exp_jn3}) and changing to polar coordinates we obtain
\begin{equation}
\EXP{J_n(r)} = n (nr^d) \th_d \int_0^2 s^{d-1} e^{-n r^d |V_d \lb \lb B_O(1) \cap B_{se_1}(1) \rb \oplus B_O(1) \rb|} \md s.
\lab{eqn:k=1_expectation_cech_up}
\end{equation}
By the Steiner's formula (\citep[(1.2)]{Schneider2008}), we have
\begin{equation}
V_d \lb \lb B_O(1) \cap B_{se_1}(1) \rb \oplus B_O(1) \rb = \th_d + \sum_{j=1}^d c_{j,d} V_j \lb B_O(1) \cap B_{se_1}(1) \rb,
\lab{eqn:Steiner}
\end{equation}
where $V_j$ are the intrinsic volumes and $c_{j,d}$ are constants depending on $j$ and the dimension $d$. For $0 \leq s \leq 2$, the lens $B_O(1) \cap B_{se_1}(1)$ contains the line segment $\ell(s)$ joining the points $\f{s}{2}e_1 - \sqrt{1 - \lb \f{s}{2} \rb^2} e_2$  and $\f{s}{2}e_1 + \sqrt{1 - \lb \f{s}{2} \rb^2} e_2$. To see this, consider the projection of the balls $B_O(1)$, $B_{se_1}(1)$ on the coordinate plane determined by the first two coordinates. Hence,
\begin{equation}
V_1 \lb B_O(1) \cap B_{se_1}(1) \rb \geq V_1(\ell(s)) = \sqrt{(2-s)(2+s)} \geq \sqrt{2(2-s)}.
\lab{eqn:lens_bound1}
\end{equation}
Using (\ref{eqn:Steiner}) in (\ref{eqn:k=1_expectation_cech_up}) and the lower bound from (\ref{eqn:lens_bound1}) we obtain
\begin{equation}
\EXP{J_n(r)} \leq n (nr^d) \th_d 2^{d-1} e^{-nr^d \th_d} \int_0^2 e^{-n r^d c_{1,d} \sqrt{2(2-s)}} \md s.
\end{equation}
Making the change of variables $u = n r^d c_{1,d} \sqrt{2(2-s)}$ we obtain
\begin{equation}
\EXP{J_n(r)} \leq C_1 \f{1}{r^d} e^{-nr^d \th_d}.
\lab{eqn:k=1_expectation_cech_up_lb2}
\end{equation}
(\ref{eqn:k=1_expectation_cech_up_lb2}) along with the fact that $\EXP{J_n(r_n^{\cC,\cU}(c_n^{\cC,\cU}))} \to e^{-\al} \in (0,\infty)$ implies that
\begin{equation} 
n(r_n^{\cC,\cU})^d \th_d + \log(r_n^{\cC,\cU})^d \leq C_2.
\lab{eqn:ub_radius_cech}
\end{equation}
Substituting for $r_n^{\cC,\cU}$ from (\ref{def:stabilising_radii}) with $k=1$ in the above inequality with $C_{\al} = \log |A^{\cC,\cU}| + \alpha -  \log m^{\cC,\cU} - \log 2$ we obtain 
\[ n(r_n^{\cC,\cU})^d \th_d + \log \lb \f{\log n + \log \log n + C_{\al}}{n c_n^{\cC,\cU}} \rb \leq C_2, \]
Adding and subtracting $\log \log n$ in the above expression, we obtain 
\[ n(r_n^{\cC,\cU})^d \th_d - \log n + \log \log n \leq C_2 - \log \lb \f{\log n + \log \log n + C_{\al}}{\log n} \rb + \log c_n^{\cC,\cU} \leq C_3 < \infty, \]
since $c_n^{\cC,\cU}$ is bounded.

To obtain the bound in the other direction, note that the lens $B_O(1) \cap B_{se_1}(1)$, $0 \leq s \leq 2$, is contained in a ball of radius $\f{\sqrt{4 - s^2}}{2}$ centered at $\f{s}{2}e_1$. Hence
\begin{equation}
V_j \lb B_O(1) \cap B_{se_1}(1) \rb \leq \lb \f{\sqrt{4 - s^2}}{2} \rb ^j V_j(B_O(1)) \leq (2-s)^{\f{j}{2}}  \; V_j(B_O(1)) .
\lab{eqn:lens_bound2}
\end{equation}
Substituting from (\ref{eqn:Steiner}) in (\ref{eqn:k=1_expectation_cech_up}) and then using the upper bound from (\ref{eqn:lens_bound2}) yields
\begin{equation}
\EXP{J_n(r)} \geq  n \th_d (nr^d) e^{n \th_d r^d} \int_0^2 s^{d-1} e^{-n r^d \sum_{j=1}^d c_{j,d} (2 - s)^{\f{j}{2}} V_j(B_O(1))} \md s.
\end{equation}
Changing variables to $u = c_{1,d} (2 - s)^{\f{1}{2}} V_1(B_O(1)) nr^d = C_1 n r^d (2 - s)^{\f{1}{2}}$ we obtain
\begin{eqnarray}
\EXP{J_n(r_n)} & \geq &  C_2 \f{n (nr_n^d)}{(nr_n^d)^2} e^{- n \th_d r_n^d} \int_0^{\sqrt{2} C_1 n r_n^d} \lb 2 - \lb \f{u}{C_1 n r_n^d} \rb^2 \rb^{d-1} u e^{- \lb u +  \sum_{j=2}^d c_{j,d}  \lb \f{u}{C_1 n r_n^d} \rb ^{j} V_j(B_O(1)) n r_n^d \rb} \md u \no \\
 & \geq & C_3 \f{1}{r_n^d} e^{- n \th_d r_n^d},
\lab{eqn:k=1_expectation_cech_up_ub2}
\end{eqnarray}
where the last inequality holds for all $n$ sufficiently large provided $n r_n^d \to \infty$ and by the dominated convergence theorem. From (\ref{eqn:k=1_expectation_cech_up_ub2}) and the fact that $\EXP{J_n(r_n^{\cC,\cU}(c_n^{\cC,\cU}))} \to e^{-\al}$ we get
\[n(r_n^{\cC,\cU})^d \th_d + \log (r_n^{\cC,\cU})^d \geq C_4. 
\]
Comparing with (\ref{eqn:ub_radius_cech}) we observe that the inequality is reversed and we have a different constant. Thus, the lower bound is obtained by a computation similar to the one following (\ref{eqn:ub_radius_cech}). 

{\it Case 3. } Finally we consider the case $p \in \cI_1$, $q = \cD$. Computations here are similar to those in the case $p = \cR$, $q = \cU$ and so we will skip some of the details. We start by observing that for $k=1$, $Q^{\cC,\cD}(O,y) = Q^{\cR,\cD}(O,y) = V_d \lb B_O(2) \cup B_y(2) \rb$, $y \in B_O(2)$ and $m^{p,\cD} =2^d\theta_d$. For $y$ of the form $s e_1$, $0 \leq s \leq 2$, we have the bounds
\begin{equation}
B_O(1) \cup B_{\lb 2 + \f{s}{2} \rb e_1 } \lb \f{s}{2} \rb \subset B_O(2) \cup B_{s e_1}(2) \subset
B_{\f{s}{2} e_1} \lb 2 + \f{s}{2} \rb .
\lab{eqn:bounds_Q_cech_down}
\end{equation}
Since $B_O(2) \cap B_{\lb 2 + \f{s}{2} \rb e_1 } \lb \f{s}{2} \rb = \emptyset$, the inclusion on the left in (\ref{eqn:bounds_Q_cech_down}) implies the following inequality.
\begin{equation}
V_d \lb B_O(2) \cup B_{s e_1}(2) \rb \geq V_d \lb B_O(2) \cup B_{\lb 2 + \f{s}{2} \rb } \lb \f{s}{2} \rb \rb = V_d \lb B_O(2) \rb + \lb \f{s}{2} \rb ^d V_d \lb B_O(1) \rb.
\lab{eqn:expr_intrinsic_vol_cech_down}
\end{equation}
Changing to polar coordinates in (\ref{eqn:exp_jn3}) and using (\ref{eqn:expr_intrinsic_vol_cech_down}) we get, 
%
%
%
\begin{equation}
\EXP{J_n(r)} \leq n (nr^d) \th_d e^{-nr^d 2^d \th_d} \int_0^2 s^{d-1} e^{-n \theta_d r^d \lb \f{s}{2} \rb ^d  } \md s.
\lab{eqn:k=1_expectation_cech_down_lb1}
\end{equation}
Making the change of variable $u = n \theta_d r^d \lb \f{s}{2} \rb ^d $ in (\ref{eqn:k=1_expectation_cech_down_lb1}) and simplifying as we did in the first two cases, we get
\begin{equation}
\EXP{J_n(r)} \leq C_1 n e^{-nr^d 2^d \th_d}
\lab{eqn:k=1_expectation_cech_down_lb2}
\end{equation}
The rest of the proof is by now a standard computation as in Part(1) (see computation following (\ref{eqn:k=1_expectation_rips_up4})). For the upper bound we use the right hand inclusion in (\ref{eqn:bounds_Q_cech_down}) to write
\begin{equation}
\EXP{J_n(r)} \geq n (nr^d) \th_d e^{-nr^d 2^d \th_d} \int_0^2 s^{d-1} e^{-n r^d \th_d \lb (2+\f{s}{2})^d - 2^d \rb} \md s.
\end{equation}
Now using the binomial expansion for $ \lb (2+\f{s}{2})^d - 2^d \rb$ and making the change of variable $u = \f{nr^d\th_ds^d}{2^d}$, we obtain
\[  \EXP{J_n(r)} \geq C_1 n e^{-nr^d 2^d \th_d} \int_0^{nr^d \th_d} e^{- u - \sum_{j=1}^{d-1} c_{j,d} \lb \f{u}{nr^d} \rb^\f{j}{d} } du. \]
Proceeding as in the proof for the lower bound \eqref{eqn:k=1_expectation_cech_up_ub2}, we derive that
\begin{equation}
\EXP{J_n(r_n)} \geq C_2 n e^{- n 2^d\th_d r_n^d},
\label{eqn:k=1_expectation_cech_down_ub2}
\end{equation}
where the last inequality holds for all $n$ sufficiently large provided $n r_n^d \to \infty$. This completes the proof of Theorem~\ref{thm:scaling_r_n_k=1}. \qed

{\bf Proof of Corollary~\ref{cor:rate_of_convergence}. } For $p \in \cI_1$ and $q \in \cI_2$, we have from Theorem~\ref{thm:scaling_r_n_k=1} that
\begin{equation}
C_1 \leq n \lb r_n^{p,q}(c_n^{p,q}) \rb ^d m^{p,q} - \log n - (1 - a^{p,q}) \log \log n \leq C_2,
\lab{eqn:bound_r_n}
\end{equation}
for some finite constants $C_1,C_2$. Substituting from (\ref{def:stabilising_radii}) in (\ref{eqn:bound_r_n}) with $C_{\al} := C_{\al}^{p,q} = \log |A^{p,q}| + \alpha -  \log m^{p,q} - \log 2$, we obtain
\[ C_1 \leq \lb \log n + \log \log n + C_{\al} \rb \f{m^{p,q}}{c_n^{p,q}} - \log n - (1 - a^{p,q}) \log \log n \leq C_2, \]
which simplifies to 
\[ C_1 \leq \lb \f{m^{p,q}}{c_n^{p,q}} - 1 \rb \log n + \lb \f{m^{p,q}}{c_n^{p,q}} - (1 - a^{p,q}) \rb \log \log n + \f{m^{p,q}}{c_n^{p,q}} C_{\al} \leq C_2. \]
The result now follows since $\f{m^{p,q}}{c_n^{p,q}} \to 1$ as $n \to \infty$ by Proposition~\ref{prop:convergence_expectation}. \qed
\subsection{Proof of result in Section \ref{sec:mon_vanishing}}
\hspace{1cm} \\ \\
{\bf Proof of Proposition \ref{prop:exp_mon_threshold} :}
	Let us fix $k \geq 1$ and $q \in \cI_2$. We first consider the case $p = \cR$. So, we shall again drop these subscripts and superscripts for rest of the calculation. Further set $\hat{J}_n(r) = J^*_n(r) - J_n(r)$. Given a $r > 0$, using \eqref{eqn:J_n^*} and Campbell-Mecke formula, we derive an upperbound for $\EXP{\hat{J}_n(r)}$. 
	\begin{eqnarray*}
	 \EXP{\hat{J}_n(r)} & = & \frac{n^{k+1}}{(k+1)!}\int_{U^{k+1}}\1[R(\x) > r]e^{-n|Q(\x,R(\x))|} \md \x  \\
	 	 & \leq & \frac{n^{k+1}}{(k+1)!}\int_{U^{k+1}}\sum_{i \neq j=1}^{k+1}\1[R(\x) > r, 2R(\x) = |x_i -x_j|]  \\
	 &  & \times e^{-n|Q(\x,|x_i-x_j|/2)|} \md \x \\
	 & \leq & \frac{n^{k+1}}{2(k-1)!}\int_{U^{k+1}} \1[2R(\x) = |x_1-x_2|]\1[|x_1-x_2| > 2r] \\ 
	 &  & \times  e^{-n|Q(\x,|x_1-x_2|/2)|} \md \x  \\
	 \mbox{($\x \to \x -(x_1,\ldots,x_1))$} & \leq & \frac{n^{k+1}}{2(k-1)!}\int_{U^k} \1[2R(O,\x) = |x_2|]\1[|x_2| > 2r] e^{-n|Q((O,\x),|x_2|/2)|} \md \x ,
	\end{eqnarray*}	
where $\x = (x_2, \ldots , x_{k+1})$ in the final expression. Changing the variable $\x \to r\x$ yields
	\begin{eqnarray*}
	 \EXP{\hat{J}_n(r)} &  \leq & \frac{n(nr^d)^{k}}{2(k-1)!} \int_{(\mR^d)^k} \1[2R(O,\x) = |x_2|,|x_2| > 2]
	 e^{-nr^d|Q((O,\x),|x_2|/2)|} \md \x.
	\end{eqnarray*}	
Changing the variable $x_2/2$ to polar co-ordinates and then $\x \to s\x$ we obtain (with $\x = (x_3, \ldots , x_{k+1})$)
	\begin{eqnarray*}
	 \EXP{\hat{J}_n(r)} 
& \leq & \frac{n2^d\theta_d(nr^d)^k}{2(k-1)!} \int_{1}^{\infty}s^{d-1} \md s \int_{B_O(2s)^{k-1}} 1[R((O,2se_1,\x)) = s] \\ 
	&  &  \times e^{-nr^d|Q((O,2se_1,\x),s))|} \md \x \\ 
& = & \frac{n2^d\theta_d(nr^d)^k}{2(k-1)!} \int_{1}^{\infty}s^{dk-1} \md s \int_{B_O(2)^{k-1}} 1[R((O,2e_1,\x)) = 1] \\
	&   & \times e^{-nr^ds^d|Q((O,2e_1,\x),1))|} \md \x.
	\end{eqnarray*}	
Since $Q(\cdot) \geq m$ for $R(\cdot) = 1$ we have
\[ \EXP{\hat{J}_n(r)}  \leq  C_1 n(n\theta_d2^dr^d)^k \int_1^{\infty}s^{dk-1} e^{-nr^ds^dm} \md s. \]
Making the change of variables $t = nr^ds^dm$ we obtain
	\begin{eqnarray*}
	 \EXP{\hat{J}_n(r)} 
& \leq & C_2 \times n\int_{nr^dm}^{\infty}t^{k-1}e^{-t}\md t \\ 
	 & \leq & C_3 \times n e^{-nr^dm}\sum_{j=0}^{k-1}\frac{(nr^dm)^j}{j!},
	\end{eqnarray*}	
	where in the last inequality we have used integral formulas for the upper gamma function. A simple substitution now yields that if $nmr_n^d = \log n + k\log \log n + w^1_n$ for some sequence $w^1_n$ bounded from below, we have that
	\[ \EXP{\hat{J}_n(r_n)} \leq C_4 e^{-w^1_n-\log \log n} \to 0. \]
   From Proposition \ref{prop:convergence_expectation_to_zero}, we know that if $nmr_n^d = \log n + k \log \log n + w^1_n$ for some sequence $w^1_n$ bounded from below, then
   \[ \EXP{J_n(r_n)} \to 0  .\]
   Thus, if  $nmr_n^d = \log n + k \log \log n + w^1_n$ for some sequence $w^1_n$ bounded from below, then
   \[ \EXP{J_{n,k}^{\cR,q,*}(r_n)} = \EXP{J_n(r_n)} + \EXP{\hat{J}_n(r_n)} \to 0, \]
   for any $q \in \cI_2$.
   
   We now consider the \Cech case i.e., $p = \cC$. In this case, the computation is a little more involved but more along the lines of that for critical points in the proof of \citep[Proposition 6.1]{Bobrowski17RSA}. Define $\hat{J}_n^{\cC}(r) = J^{\cC,q,*}_n(r) - J^{\cC,q}_n(r)$. By definition of the \v{C}ech complex, $R^{\cC}(x_0,\ldots,x_k) = \inf \{r : \cap_{i=0}^kB_r(x_i) \neq \emptyset\}$ and further we have that  $\{C(x_0,\ldots,x_k)\} = \cap_{i=0}^kB_{R^{\cC}(x_0,\ldots,x_k)}(x_i)$ for some point $C(x_0,\ldots,x_k) \in \mR^d$. Now using this observation and proceeding as in the Vietoris-Rips complex case using translation and scaling we have that (dropping the superscripts $\cC,q$ as usual)
\begin{equation*} 
\EXP{\hat{J}_n(r)} \leq  \frac{n^{k+1}r^{dk}}{(k+1)!}\int_{(\mR^d)^k}\1[R(O,\x) > 1]e^{-nr^d|Q((O,\x),R(O,\x))|} \md \x 
\end{equation*}
The RHS in the above equation is exactly of the form \citep[(8.8)]{Bobrowski17RSA} with $h_1(O,\x)$ there replaced by $\1[R(O,\x) > 1]$ and $\theta_d R(O,\x)^d$ replaced by $|Q((O,\x),R(O,\x))|$ in the exponent. Observe that both $|Q((O,\x),R(O,\x))|$  and $R(O,\x)$ are rotation invariant and also $|Q((O,s\x),R(O,s\x))| = s^d|Q((O,\x),R(O,\x))|$ for any $s > 0$. So, we can now follow the derivations in \citep[(8.8)-(8.10)]{Bobrowski17RSA} and using the bound that $Q((O,\x),1) \geq m$ derive that 
$$ \EXP{\hat{J}_n(r)} \leq C_1 n (nr^d)^k \int_1^{\infty}s^{dk-1}e^{-nr^ds^dm} \md s.$$
The above integral can be simplified and evaluated as in the Vietoris-Rips case above to obtain that
$$ \EXP{\hat{J}_n(r)} \leq C_2 \times ne^{-nr^dm}\sum_{j=0}^{k-1}\frac{(nr^dm)^j}{j!}.$$
Thus again combining with the Proposition \ref{prop:convergence_expectation_to_zero}, we have that if $nmr_n^d = \log n + k \log \log n + w^1_n$ for some sequence $w^1_n$ bounded from below, then $\EXP{J_{n,k}^{\cC,q,*}(r_n)} \to 0$ for any $q \in \cI_2$. \qed
\subsection{Proofs of results in Section \ref{sec:preview}}
\hspace{1cm} \\ \\
{\bf Proof of Theorem~\ref{thm:coarse_phase_transition}. } Fix $p \in \cI_1, q \in \cI_2$ and we shall drop the superscripts $p,q$ in the rest of the proof. Let $nm_kr_n^d = (1+\ep) \log n$. 
Substituting in (\ref{eqn:exp_jn3}) and observing that $|Q(O,y)| \geq m_k$ and $A$ is bounded, we obtain
\begin{eqnarray}
\EXP{J_n(r_n)} & \leq & C_1 n(nr_n^d)^k e^{- nr_n^d m_k} \no \\
 & \leq & C_2 \f{(\log n)^k}{n^{\ep}}. 
\lab{eqn:bound_exp}
\end{eqnarray}
By the Markov's inequality and the bound obtained in (\ref{eqn:bound_exp}), we have
\[ \pr{J_n(r_n) \geq 1} \leq C_2 \f{(\log n)^k}{n^{\ep}} \to 0, \]
as $n \to \infty$. This proves the second assertion in \eqref{eqn:coarse_phase_transition}.

Let $nm_k^{p,q}r_n^d = (1-\ep) \log n$. To prove the first assertion we use the second moment approach. Since  
\[ \pr{J_n(r_n) \geq 1} \geq \f{\lb \EXP{J_n(r_n)} \rb^2}{\EXP{J_n(r_n)^2}}, \]
to prove the first assertion it suffices to show that

\begin{equation}
\f{\lb \EXP{J_n(r_n)} \rb^2}{\EXP{J_n(r_n)^2}} \to 1, \qquad \mbox{as } n \to \infty. 
\lab{eqn:second_moment}
\end{equation}
To this end, we evaluate $\EXP{J_n^2}$. From (\ref{def:J}) we can write 
\begin{eqnarray}
J_n(r_n)^2 & = & C_1 \sum_{\x, \y \in \cP^{k+1}_n} h(\x,r_n) h(\y,r_n) \1 \left[\cP_n \lb Q(\x,r_n) \cup Q(\y,r_n) \rb = 0 \right] \no \\
 & = & \sum_{j = 0}^{k+1} J_n^{(j)},
\lab{eqn:split_J_n_square}
\end{eqnarray}
where
\begin{equation}
J_n^{(j)} =  C_1 \sum_{\stackrel{\x, \y \in \cP^{k+1}_n,}{|\x \cap \y| = j}} h(\x,r_n) h(\y,r_n) \1 \left[\cP_n \lb Q(\x,r_n) \cup Q(\y,r_n) \rb = 0 \right], \qquad j = 0,1,2, \ldots ,(k+1),
\end{equation}
is the contribution to $J_n^2$ when the two complexes share $j$ vertices. For $j=1, \ldots , (k+1)$, we have by the Campbell-Mecke formula
\begin{eqnarray}
\EXP{J_n^{(j)}} & \leq & C_2 n^{2k+2 - j} \int_{\x \in U^{(k+1)}} \int_{\z \in U^{(k+1-j)}} h(\x,r_n) h(\y,r_n) e^{- |Q(\x,r_n) \cup Q(\y,r_n)|} \md \x \md \z \no \\
 & \leq & C_3 (n r_n^d)^{k+1-j} n^{k+1} \int_{\x \in U^{(k+1)}} h(\x,r_n) e^{- |Q(\x,r_n)|} \md \x,
\lab{eqn:bound_jnj}
\end{eqnarray}
where $\y = (x_1, \ldots ,x_j,\z), \z = (z_1,\ldots,z_{k+1-j})$ and the last inequality in (\ref{eqn:bound_jnj}) follows from the restriction that the variables $z_i$, $i=1, \ldots , (k+1-j)$ all lie within a ball of radius $6 k r_n$ from $x_1$ and the integrand in the $\z-$variables is bounded by one. Comparing the right hand side of the last expression in (\ref{eqn:bound_jnj}) with (\ref{eqn:exp_jn1}) and using the definition of $r_n$ we obtain
\begin{equation}
\EXP{J_n^{(j)}} \leq C_4 \lb \log n \rb ^{k+1-j} \EXP{J_n}, \qquad j=1,2, \ldots (k+1).
\lab{eqn:bound_jnj1}
\end{equation}
Now consider $\EXP{J_n^{(0)}}$. By the Campbell-Mecke formula we have
\begin{equation}
\EXP{J_n^{(0)}}  = [(k+1)!]^{-2} n^{2k+2} \int_{\x \in U^{(k+1)}} \int_{\z \in U^{(k+1)}} h(\x,r_n) h(\z,r_n) e^{- |Q(\x,r_n) \cup Q(\z,r_n)|} \md \x \md \z.
\lab{eqn:jn0}
\end{equation}
Divide the inner integral in (\ref{eqn:jn0}) into two parts, one over the region where $\min_{1 \leq r,s \leq (k+1)} |x_r - z_s| \leq 6 k r_n$ and the second its complement. Over the first region we proceed as in (\ref{eqn:bound_jnj}), (\ref{eqn:bound_jnj1}) to obtain the bound 
\[ n^{2k+2} \int_{\x \in U^{(k+1)}} \int_{\z \in U^{(k+1)}} \1[\min_{1 \leq r,s \leq (k+1)} |x_r - z_s| \leq 6 k r_n] \;
h(\x,r_n) h(\z,r_n) e^{- |Q(\x,r_n)|} \md \x \md \z \]
\begin{equation} \leq C_5 \lb \log n \rb^{k+1} \EXP{J_n}.
\lab{eqn:jn0_bd1}
\end{equation}
Over the region where $\min_{1 \leq r,s \leq (k+1)} |x_r - z_s| > 6 k r_n$, we have $|Q(\x,r_n) \cup Q(\z,r_n)| = |Q(\x,r_n)| + |Q(\z,r_n)|$ which yields the bound
\[ [(k+1)!]^{-2} n^{2k+2} \int_{\x \in U^{(k+1)}} \int_{\z \in U^{(k+1)}} \1[\min_{1 \leq r,s \leq (k+1)} |x_r - z_s| > 6 k r_n] \;
h(\x,r_n) h(\z,r_n) e^{- |Q(\x,r_n)| + |Q(\z,r_n)|} \md \x \md \z \]
\begin{equation} \leq \lb \EXP{J_n} \rb ^2.
\lab{eqn:jn0_bd2}
\end{equation}
From (\ref{eqn:split_J_n_square}), (\ref{eqn:bound_jnj1}) - (\ref{eqn:jn0_bd2}) we obtain

\begin{equation}
\EXP{J_n^2} \leq C_7 \lb \log n \rb^{k+1} \EXP{J_n} + \lb \EXP{J_n} \rb ^2.
\lab{eqn_exp_jn_square_bound}
\end{equation}
Choose $\del > 0$ sufficiently small so that $\f{(m_k + \del)(1-\ep)}{m_k} = 1 - \f{\ep}{2}$. With this choice of $\del$, substituting for $r_n$ in (\ref{eqn:exp_jn3}) we obtain
\begin{eqnarray}
\EXP{J_n} & \geq & \frac{n(n r_n^d)^{k}}{(k+1)!}\int_A \1[|Q(O,\y)| \leq m_k + \del] \; e^{- n r_n^d |Q(O,\y)|} \md \y \no \\
 & \geq & C_8 n (\log n)^k e^{- \f{(m_k + \del)(1-\ep) \log n}{m_k}} \no \\
 & = & C_8 n^{\f{\ep}{2}} (\log n)^k.
\lab{eqn:lb_exp_jn}
\end{eqnarray}
It now follows from (\ref{eqn_exp_jn_square_bound}) and (\ref{eqn:lb_exp_jn}) that
\[ \liminf_{n \to \infty} \f{\lb \EXP{J_n(r_n)} \rb^2}{\EXP{J_n(r_n)^2}} \geq 1. \]
This proves \eqref{eqn:second_moment} and hence the first assertion in \eqref{eqn:coarse_phase_transition}. For any $m \geq 1$, we have by Chebyshev's inequality
\[  P(J_n(r_n) \leq m) \leq \f{\VAR{J_n(r_n)}}{(\EXP{J_n(r_n)} - m)^2}  \to 0,  \]
where the convergence follows from \eqref{eqn:lb_exp_jn} and \eqref{eqn:second_moment}. This proves that $J_n(r_n) \stackrel{P}{\to} \infty$.  \qed

{\bf Proof of Theorem \ref{thm:coarse_phase_transition1}:} The first statement in \eqref{eqn:coarse_phase_transition1} follows trivially from the corresponding statement in \eqref{eqn:coarse_phase_transition} and the second statement now follows from Proposition \ref{prop:exp_mon_threshold} and Markov's inequality. \qed

 {\bf Proof of Theorem \ref{thm:finite_components}:}
   
   Fix $p \in \cI_1, q \in \cI_2, k \geq 1$ and $L \geq 1$. We shall now onwards drop  superscripts $p,q$ in the rest of the proof except to avoid ambiguity. For $M \geq 1$ let $\Delta(L,M)$ denote the set of feasible (up/down)-connected graphs $\Gamma$ formed by $L$ $k$-faces of the \v{C}ech or Vietoris-Rips complex such that there are a total of $M$ vertices in the $k$-faces and each of the $M$ vertices is present in at least one of the $k$-faces. More precisely, $\Gamma \in \Delta(L,M)$ if there exists $\{x_1,\ldots,x_M\} \subset \mR^d$ with $S_k(\{x_1,\ldots,x_M\},1) = L, G_k(\{x_1,\ldots,x_M\},1) \cong \Gamma$ and further each $x_i, 1 \leq i \leq M$ belongs to at least one $k$-face. Note that it is possible that $\Delta(L,M) = \emptyset$ for certain choices of $M$ and $L$ either due to the combinatorics or the geometry. Trivially, $\Delta(L,M) = \emptyset$ for $M > L(k+1)$ and $M \leq k$. Hence setting $\Delta(L) = \cup_{M=k+1}^{L(k+1)}\Delta(L,M)$, we see that $\Delta(L)$ is the set of all feasible (up/down)-connected graphs that can be formed on $L$ faces. Note that both $\Delta(L,M)$ and $\Delta(L)$ depend on $p$ and $q$ but we omit the same. Thus, we have that 
\begin{equation}
\label{eqn:jnkrl}
J_{n,k}(r,L) = \sum_{M=k+1}^{L(k+1)} \sum_{\Gamma \in \Delta(L,M)} \bar{J}_{n,k}(r,\Gamma),
\end{equation}
   where $\bar{J}_{n,k}(r,\Gamma)$ is the number of induced $\Gamma$ components in $G_k(\cP_n,r)$ formed by $M$ vertices i.e.,
 \begin{equation}
\label{eqn:jnkrgamma}
\bar{J}_{n,k}(r,\Gamma) := \sum_{\{X_1,\ldots,X_M\} \subset \cP_n}\1[G_k(\{X_1,\ldots,X_M\},r) \cong \Gamma]\1[\cP_n(\cup_{i=1}^{L}Q((\X^i),r)) = \{X_1,\ldots,X_M\}],\
\end{equation}
   where $\X^i$ is a $(k+1)$-subset of $\{X_1,\ldots,X_M\}$ such that the $\X^i, i = 1,\ldots,L$ are the vertices in $G_k(\{X_1,\ldots,X_M\},r)$ i.e., the $k$-faces in the corresponding geometric complex. Since $\Delta(L)$ is a finite set, it is enough if we show that for all $\Gamma \in \Delta(L,M)$,
   $$\EXP{\bar{J}_{n,k}(r_n,\Gamma)} \to 0$$
   for $r_n$ such that $nm_kr_n^d = (1 + \epsilon) \log n$ for any $\epsilon > 0$. 

   Set $h_{\Gamma,r} := \1[G_k(\{x_1,\ldots,x_M\},r) \cong \Gamma]$ and $h_{\Gamma} := h_{\Gamma,1}$. Since $\Gamma$ is connected, we note that there exists a $K > 0$ (possibly depending on $M,L$) such that $h_{\Gamma}(O,x_2,\ldots,x_M) =  0$ if $\max_{i=2,\ldots,M}|x_i| > K$. 
   
   Further, whenever $G_k(\{x_1,\ldots,x_M\},r) \cong \Gamma$, we denote the $L$ vertices (i.e., $k$-faces) by $\x^1,\ldots,\x^L$. Let $r > 0$. As usual, we start with the Campbell-Mecke formula and then use translation and scaling relations in the below derivation :
   \begin{eqnarray*}
   \EXP{\bar{J}_{n,k}(r,\Gamma)} & = & \frac{n^M}{M!}\int_{U^M}h_{\Gamma,r}(x_1,\ldots,x_M)e^{-n|\cup_{i=0}^LQ(\x^i,r)|}\md x_1 \ldots \md x_M \\
   (\mbox{change $x_i \to x_i + rx_1, i \geq 1$}) & \leq & \frac{n(nr^d)^{M-1}}{M!}\int_{(\mR^d)^{M-1}}h_{\Gamma}(O,x_2,\ldots,x_m)e^{-nr^d|\cup_{i=0}^LQ(\x^i)|} \md x_2 \ldots \md x_M \\
   (\mbox{by $|\cup_{i=0}^LQ(\x^i)| \geq m_k$}) & \leq & \frac{n(n\theta_dK^dr^d)^{M-1}}{M!} e^{-nr^d m_k}.
   \end{eqnarray*}

  Now choosing $r_n$ such that $nm_kr_n^d = (1 + \epsilon) \log n$ for an $\epsilon > 0$, we have using the above bound that
  $$ \EXP{\bar{J}_{n,k}(r_n,\Gamma)} \leq  (\frac{\theta_dK^d(1+\epsilon)^d}{m_k})^{M-1}n^{-\epsilon}(\log n)^{M-1}  \to 0.$$ 
\qed 

{\bf Proof of Proposition~\ref{prop:exp_Jn1}. } The results are a straightforward consequence of the inequalities (\ref{eqn:k=1_expectation_rips_up4}), \eqref{eqn:k=1_expectation_rips_up5},      (\ref{eqn:k=1_expectation_cech_down_lb2}), (\ref{eqn:k=1_expectation_cech_down_ub2}), (\ref{eqn:k=1_expectation_cech_up_lb2}) and (\ref{eqn:k=1_expectation_cech_up_ub2}) obtained in the proof of Theorem~\ref{thm:scaling_r_n_k=1}. \qed

{\bf Proof of Proposition \ref{prop:exp_Jn1s}:}

Fix a $p \in \cI_1, q \in \cI_2$. Recall that $\hat{J}_{n,1}(r) = J^*_{n,1}(r) - J_{n,1}(r)$. We note that $R(x_1,x_2) = |x_1-x_2|/2$. Again, we know asymptotics of $\EXP{J_{n,1}(r)}$ from Corollary \ref{cor:rate_of_convergence} and Theorem \ref{thm:scaling_r_n_k=1}. So, we shall only derive asymptotics for $\EXP{\hat{J}_{n,1}(r)}$. Again, starting with Campbell-Mecke formula and using translation, change to polar coordinates as in the above calculations
   \begin{eqnarray*}
    \EXP{\hat{J}_{n,1}(r)} & = & \frac{n^2}{2}\int_{U^2}\1[R(\x) > r]e^{-n|Q(\x,R(\x))|} \md \x \\
     & = & \frac{n^2}{2}\int_{U} \md x_0 \int_{U-x_0}\1[|x-x_0| > 2r]e^{-n|Q((x_0,x),|x|/2)|} \md x \\
   ((x_0,x) \to (O,x-x_0))  & = & \frac{n^2}{2} \int_{U}\1[|x| > 2r]e^{-n|Q((O,x),|x|/2|)} \md x \\
   (\mbox{change $x/2r$ to polar co-ordinates})      & =  & n^2\theta_d2^{d-1}r^d \int_{1}^{\infty} s^{d-1} e^{-nr^d|Q((O,2se_1),s)|} \md s \\
     & =  & n^2\theta_d2^{d-1}r^d \int_1^{\infty} s^{d-1} e^{-nr^ds^d |Q((O,2e_1),1)|} \md s \\
(\mbox{by definition of $m_2$}) & =  & n^2\theta_d2^{d-1}r^d \int_1^{\infty} s^{d-1} e^{-nr^ds^dm_2} \md s \\
     & = &   \frac{2^{d-1}\theta_d}{dm_2} ne^{-nr^dm_2}.
   \end{eqnarray*}
   Thus, combining with Proposition \ref{prop:exp_Jn1}, the proof is complete.  \qed 

\subsection{Proofs of results in Section \ref{sec:poisson_convg}} 

Theorem~\ref{thm:rips-Poisson-cgs} is proved using the criterion derived in \citep[Theorem 3.1]{Penrose15}, a simpler version of which is stated below. In order to state this Theorem, we need some notation. Let $\eta$ be a finite Poisson point process in $\mR^d$ with intensity measure $\mu$ and $\cN$ be the space of all finite subsets of $\mR^d$ with the sigma-algebra on $\cN$ generated by the functions $\xi \to |\xi \cap B|$ for all bounded Borel sets $B \subset \mR^d$. Let $k \in \mN$ and let $f:(\mR^d)^k \times \cN \to \{0,1\}$ be a measurable function. For any $\xi \in \cN$, set
\[ F(\xi) := \sum_{\psi \subset \xi: | \psi | = k} f(\psi, \xi\setminus \psi). \]
For $x_1, \ldots ,x_k \in \mR^d$ set $p(x_1, \ldots , x_k) = \EXP{f(\{ x_1, \ldots ,x_k \}, \eta)}$. 

\begin{theorem}(\citep[Theorem 3.1]{Penrose15})
\label{thm:Penrose15}	
Let $W = F(\eta)$ with $\eta$ and $F$ as defined above. Suppose that $w: (\mR^d)^k \to [0, \infty)$ is a measurable function, and that for $\mu^k$-almost every $\x = (x_1, \ldots ,x_k) \in (\mR^d)^k$ with $p(\x) > 0$ we can find coupled random variables $U^{\x}, V^{\x}$ such that 
\begin{itemize}
\item $U^{\x} \stackrel{d}{=} W$
\item $1 + V^{\x} \stackrel{d}{=} F(\cup_{i=1}^k \{x_i \} \cup \eta) \bigg| f(\{x_1, \ldots ,x_k\}, \eta) = 1$
\item $\EXP{|U^{\x} - V^{\x}|} \leq w(\x).$
\end{itemize}
Then the total variation distance between the law of $W$ and a Poisson random variable with mean $\EXP{W}$ satisfies
\[ d_{TV}(W, Poi(\EXP{W})) \leq \f{1 \wedge (\EXP{W})^{-1}}{k!} \int w(\x)p(\x) \mu^k(\md\x). \]
\end{theorem}
The following geometric lemma is crucial in the proofs of Theorems \ref{thm:rips-Poisson-cgs} and  \ref{thm:rips-Poisson-cgs_finite}. The lack of such a geometric lemma hinders extending these results to the \v{C}ech complex. 
\begin{lemma} 
\label{lem:geometric} For $\delta \geq 0, j \leq k$, define
 $$ D_{j,\delta} :=   \{ (\x,\z) \in B_O(2)^{k} \times B_O(6)^{k-j+1}: 2- \delta \leq |x_i|, |x_i - x_j| , |z_i|, |x_i - z_j|, |z_i - z_j|, \forall i \neq j, h((O,\x)) = h(\y) = 1 \},$$ 
where  $\x = (x_2, \ldots , x_{k+1}) \in \mR^{dk}$, $\z = (z_1,\ldots z_{k-j+1}) \in \mR^{d(k-j+1)}$ and \\ $\y = (x_{k-j+2}, \ldots ,x_{k+1}, z_1,\ldots z_{k-j+1}) \in \mR^{d(k+1)}$. Then there exists a $\delta_0 >0$ such that for any $0 \leq\delta < \delta_0$ we can find a $\beta := \beta(\delta) > 0$ for which $|Q(\y) \setminus Q(O,\x)| \geq \beta$ on the set $D_{j,\delta}$.
\end{lemma}
{\bf Proof.} Since the function $|Q(\y) \setminus Q(O,\x)|$ is continuous in $(\x,\z)$, it suffices to show the result with $\delta = 0$. Then any $(\x,\z)  \in  D_{j,0}$ must satisfy the following conditions. Firstly, since $ h((O,\x)) = h(\y) = 1$, we have that $|x_i|, |x_i - x_{\ell}| , |z_i - z_{\ell}| = 2, \forall i \neq {\ell}$ and also $|x_i - z_{\ell}| = 2, \forall {\ell}$ and $\forall i \in \{k-j+2,\ldots,k+1\}$. Secondly, $|z_{\ell}|, |x_i - z_{\ell}| \geq 2, \forall {\ell}$ and $\forall i \in \{2,\ldots,k-j+1\}.$ Since $Q(O,\x) \subset B_O(2)$, it suffices to show $|Q(\y) \setminus B_O(2)| \geq \beta$ for some $\beta > 0$.  

We now state two claims which will be proven later.
\begin{enumerate}
\item[{\em Claim 1:}] Define $\delta_1 := \min\{ \max_{i \neq {\ell}} |x_i - x_{\ell}| :  x_1,\ldots,x_{d+2} \in \mR^d, |x_i - x_{\ell}| \geq 2, \, \, \forall i \neq {\ell} \} - 2$. The first claim is that $\delta_1 > 0$. 
\item[{\em Claim 2:}] Let $\delta_1 >0$ be as in Claim 1. If $x_1,\ldots,x_{d+1} \in \mR^d$ are such that  $|x_i - x_{\ell}| = 2, \, \, \forall i \neq {\ell}$, then there exists an $x \in \mR^d$ such that $|x - x_1| \geq 2 + \delta_1$ and $|x - x_i| = 2$ for all $i \in \{2,\ldots,d+1\}$. 
\end{enumerate}
Using the above two claims, we now complete the proof.

First consider the case when $k = d$. Since $\delta = 0$, $x_{k-j+2},\ldots, x_{k+1}, z_1 , \ldots z_{k-j+1}$ are at distance exactly  two from each other.  Hence, from Claim 1, we have that $|z_i| > 2 +\delta_1$ for some $i \in \{1,\ldots k-j+1\}$. Without loss of generality, let us assume that $|z_1| > 2 +\delta_1$. Hence $B_{z_1}(\delta_1/2) \cap B_O(2) = \emptyset$ and also since $z_1 \in Q(\y)$, we have by convexity that $|B_{z_1}(\delta_1/2) \cap Q(\y)| \geq \beta > 0$ for some $\beta > 0$. Thus, we get that $|Q(\y) \setminus B_O(2)| \geq \beta$ for some $\beta > 0$.

Next, consider the case when $k < d$. Take the points $x_{k-j+2}, \ldots , x_{k+1}, z_1 , \ldots z_{k-j+1}$.  Since we are interested in minimizing  $|Q(\y) \setminus B_O(2)|$, we can assume that $|z_i| = 2$ for all $i \in \{1,\ldots k+1\}$. Further, if $k < d-1$ choose additional points $\ze_1, \ldots \ze_{d-k-1}$ all on the boundary of $B_O(2)$ so that $(O, \bze)$ forms a $d$-simplex with side lengths $2$ where $\bze = (x_{k-j+2}, \ldots , x_{k+1}, z_1 , \ldots z_{k-j+1}, \ze_1, \ldots , \ze_{d-k-1})$. When $k = d-1$, we set  $\bze = (x_{k-j+2}, \ldots , x_{k+1}, z_1 , \ldots z_{k-j+1})$ and observe that it still holds that $(O, \bze)$ forms a $d$-simplex with side lengths $2$. From Claim 2, we can choose $x \in \mR^d$ such that $|x| \geq 2 + \delta_1$ and $(\bze,x)$ forms a $d$-simplex with side lengths $2$. Since $x \in Q(\y)$, we can argue as in the case $k =d$ by convexity that  $|Q(\y) \setminus B_O(2)| \geq |B_{x}(\delta_1/2) \cap Q(\y)| \geq \beta$ for some $\beta > 0$.  

This completes the proof except the two claims which will be proven next. 

{\em Proof of Claim 1:} Without loss of generality, we can choose $x_1 = O$. We can further assume that $|x_i| \leq 3$ for all $i = 2,\ldots d+1$ as minimum will be attained by such a configuration of points. Since $\max_{i \neq {\ell}} |x_i - x_{\ell}|$ is a continuous function of $x_2,\ldots x_{d+1}$ on $B_O(3)^d$, the minimum $\delta_1$ is attained. If $\delta_1 = 0$, we have a contradiction that there is a configuration of $(d+2)$ points which form a $(d+1)$-simplex in $\mR^d$ with side-lengths $2$.

{\em Proof of Claim 2:} To show this we will make a specific choice of $\x = (x_1,\ldots x_{d+1})$ and show existence of $x$. Any other choice will be a rotation and translation of this configuration. To simplify notation we relabel $\x$ to be $(O,x_1, \ldots,x_d)$. Let $e_i$ $i=1, \dots ,d$ be the unit vectors along the coordinate axes. Let $\bv = \sum_{i=1}^d e_i$ and write for $i = 1, \ldots , d$, $x_i = a e_i + b \bv$. Since $\x$ forms a simplex with side lengths two, the constants $a, b$ must satisfy the following two conditions: $(i)$ $|x_i| = 2$, $i = 1, \ldots , d$ and $(ii)$ $|x_i - x_{\ell}| = 2$, $i,\ell = 1,2, \ldots ,d$, $i \ne \ell$. In other words, to obtain $x_i$ we start with the vector of length two along $e_i$ and rotate it towards $\bv$.

From condition $(ii)$ above, it follows that $a = \sqrt{2}$ and from $(i)$ $b$ is the positive solution of the quadratic equation $(a+b)^2 + (d-1) b^2 = 4$ and thus 
\begin{equation}
b = \f{- \sqrt{2} + \sqrt{2} \sqrt{1+d}}{d} \qquad \Rightarrow \qquad \f{a + db}{\sqrt{d}} = \sqrt{2} \sqrt{\f{1 + d}{d}} > 1.
\label{eqn:exp_b}
\end{equation}
Denote by $C = \f{\sum_{i=1}^d x_i}{d} = \lb \f{a}{d} + b \rb \bv$ the centroid of the points in $\bze$. Choose $x = 2C$. Thus $C$ is on the hyperplane containing the points in $\bze$ and $O,C$ and $x$ are collinear with $C$ being the mid point of the line segment joining $O$ and $x$. From (\ref{eqn:exp_b}) we obtain
\[ |x| = 2  \lb \f{a}{d} + b \rb \sqrt{d} = 2 \f{a + db}{\sqrt{d}} > 2. \]
$O$ is at a distance two from all the points in $\x$ and by symmetry, the point $x$ is also at a distance two from all the points in $\bze$. Thus,  $(\x,x)$ satisfy the assumptions of Claim 1  and so $|x| > 2 + \delta_1$ by Claim 1. \qed

{\bf Proof of Theorem~\ref{thm:rips-Poisson-cgs}.}
Throughout this proof we take $p=\cR$ and $q = \cU$ and so omitting the superscripts, we will denote $r_n^{\cR,\cU}(c_n)$ by $r_n$ and $J^{\cR,\cU}_n(r_n^{\cR,\cU}(c_n))$ by $J_n$. Recall that the sequence $\{c_n\}$ satisfies
\begin{equation}
\beta_n = \EXP{J_n(r_n(c_n))} \to e^{-\alpha} \qquad \mbox{as } n \to \infty.
\lab{eqn:def_beta_n}
\end{equation}

Recall that $Q^{\cR,\cU}(\x,r_n) = \cap_{i=1}^{k+1} B_{x_i}(2r_n)$ and let $m = m^{\cR,\cU}$. Also, we set $h_n(\y) = h(\y,r_n), Q_n(\y) = Q(\y,r_n), \tilde{Q}_n(\y) = Q_n(\y) \setminus \{\y\}$. The proof follows by verifying the criterion given in Theorem~3.1. To invoke this criterion, take $\eta = \cP_n, f(\y,\eta) = h_n(\y)1[\eta(Q_n(\y)) = 0].$ So, $W_n = F(\cP_n)$ is the number of isolated Rips $k$-complexes in the graph $G_k(\cP_n,r_n)$ and more explicitly,
\begin{equation}
W_n = J_n(r_n(c_n))  = F(\cP_n) = \sum_{\y \in \cP_n^{(k+1)}}f(\y,\cP_n \setminus \y) = \sum_{\y \in \cP_n^{(k+1)}}h_n(\y)1[\cP_n(\tilde{Q}_n(\y)) = 0].
\lab{eqn:def-W}
\end{equation}
For any $\x \in (\mR^d)^{k+1}$, set $\cP_n^{\x} = (\cP_n \cap Q_n(\x)^c) \cup \{\x\}$. Set $U_n^{\x} = W_n$ and define $V_n^{\x}$ as
$$ V_n^{\x} = \sum_{\stackrel{\y \in (\cP_n^{\x})^{(k+1)}}{\y \ne \x}}f(\y,\cP_n^{\x}) = \sum_{\stackrel{\y \in (\cP_n^{\x})^{(k+1)}}{\y \ne \x}}h_n(\y) 1[\cP_n^{\x}(\tilde{Q}_n(\y)) = 0],$$
where $\y \neq \x$ denotes that $\y$ differs from $\x$ in at least one co-ordinate.
Let $\x$ be such that $p_n(\x) := \EXP{f(\x,\cP_n)} > 0.$ In particular, this implies $h_n(\x) = 1$.
$$1 + V_n^{\x} = f(\x,\cP_n^{\x} \setminus \x) + V_n^{\x} = F(\cP_n^{\x}) \stackrel{d}{=} F(\cP_n \cup \x)\big|\{f(\x,\cP_n) = 1\}.$$
The first equality follows because $h_n(\x) = 1$ and $\cP_n^{\x}(\tilde{Q}_n(\x)) = 0$, the second equality follows from definition of $F(\cP_n^{\x})$ (see \eqref{eqn:def-W}) and the third equality follows because $\cP_n \cup \{\x\} \big| \{f(\x,\cP_n) = 1\} = \cP_n \cup \{\x\} \big| \{ \cP_n(\tilde{Q}_n(\x)) = 0 \} \stackrel{d}{=}  \cP_n^{\x}$. We define
\[W_n^{(1)}(\x)  :=  \sum_{\stackrel{\y \in \cP_n^{(k+1)}}{Q_n(\y) \cap Q_n(\x) \neq \emptyset}}f(\y,\cP_n), \qquad
W_n^{(2)}(\x)  :=   \sum_{\stackrel{\y \in (\cP_n^{\x})^{(k+1)}, \y \ne \x}{Q_n(\y) \cap Q_n(\x) \neq \emptyset}}f(\y,\cP_n^{\x}), \]
\begin{equation}
W_n^{(3)}(\x)  :=  \sum_{\stackrel{\y \in \cP_n^{(k+1)}}{Q_n(\y) \cap Q_n(\x) = \emptyset}}f(\y,\cP_n) = \sum_{\stackrel{\y \in (\cP_n^{\x})^{(k+1)}}{Q_n(\y) \cap Q_n(\x) = \emptyset}}f(\y,\cP_n^{\x}), 
\end{equation}
where the last equality follows by observing that $\y \neq \x$ and $f(\y,\cP_n^{\x}) = f(\y, \cP_n)$ if $Q_n(\y) \cap Q_n(\x) = \emptyset$. Now, we can write $U_n^{\x} = W_n^{(1)}(\x) + W_n^{(3)}(\x)$ and $V_n^{\x} = W_n^{(2)}(\x) + W_n^{(3)}(\x).$ This yields $|U(\x) - V(\x)| \leq W_1(\x) + W_2(\x)$. 

We let $w_n(\x) = w_n^{(1)}(\x) + w_n^{(2)}(\x)$ with $w_n^{(i)}(\x) = \EXP{W_n^{(i)}(\x)}$, $i=1,2$. Then applying Theorem~3.1 we obtain
\begin{equation}
d_{TV}(J_n, Poi(\beta_n)) \leq  \frac{1 \wedge \beta_n^{-1}}{(k+1)!}(I_1 + I_2),
\lab{eqn:bound_dTV}
\end{equation}
where
\begin{equation}
I_i = n^{k+1} \int_{U^{k+1}} w_n^{(i)}(\x)p_n(\x) \md \x, \qquad i=1,2.
\end{equation}
The result now follows from (\ref{eqn:def_beta_n}) and (\ref{eqn:bound_dTV}) provided we show that $I_i \to 0$ as $n \to \infty$ for $i=1,2$. Recall that 
\[   p_n(\x) = \EXP{f(\x,\cP_n)} = h_n(\x)e^{-n |Q_n(\x)|}.\]
By the Campbell-Mecke formula applied to $w_1(\x)$ as in (\ref{eqn:exp_jn1}), (\ref{eqn:exp_jn2}) and noting that $Q_n(\y) \cap Q_n(\x) \ne \emptyset$ as well as setting $Q(O,y) := Q((O,\y),1)$, we obtain
\begin{eqnarray}
I_1 & = & C_1 n^{k+1} \int_{U^{k+1}}\md \x \, \, p_n(\x) \, \, n^{k+1} \int_{\{\y: Q_n(\y) \cap Q_n(\x) \ne \emptyset\}} p_n(\y) \md \y \no \\
 & \leq & C_1 n^{2(k+1)} \int_{\x \in U^{k+1}} \int_{\y \in B_{x_1}(6r_n)\times U^k} h_n(\x) h_n(\y) e^{-n (|Q_n(\x)| + |Q_n(\y)|)} \md \x \md \y \no \\ 
 & \leq & C_2 r_n^d \left(n(nr_n^d)^{k} \int_{(r_n^{-1}U)^k} h(O,\y) e^{-n r_n^d |Q(O,\y)|}  \md \y \right)^2 = C_2 r_n^d \beta_n^2 \to 0, \label{eqn:cgs_I_1}
\end{eqnarray}
from (\ref{eqn:exp_jn3}), (\ref{eqn:def_beta_n}) and the fact that $r_n \to 0$ as $n \to \infty$. 

To analyse $I_2$ we will write it as a sum depending on the number of coordinates common to $\x$ and $\y$.
\[ I_2 = \sum_{j=0}^{k} I_{2j}, \]
where
\begin{equation} 
I_{2j} = \binom{k+1}{j} n^{2(k+1)-j} \int_{\x \in U^{k+1}} \int_{\z \in (U \setminus Q(\x,r_n))^{k-j+1}} h_n(\x) h_n(\y) e^{-n |Q_n(\x) \cup Q_n(\y)|} \md \x \md \z,
\lab{eqn:I_2j}
\end{equation}
$\y = (x_{k-j+2}, \ldots ,x_{k+1}, z_1,\ldots z_{k-j+1})$ and $\z = (z_1,\ldots z_{k-j+1})$. Note that $\x, \y$ have $j$ coordinates $(x_{k-j+2}, \ldots ,x_{k+1})$ in common. Since the metric is toroidal, for any $x_1 \in U$ the integration with respect to the remaining variables yields a function that does not depend on $x_1$. Hence we can fix the first variable to be the origin $O$. Bounding (\ref{eqn:I_2j}) as in (\ref{eqn:cgs_I_1}), we obtain
\begin{equation}
I_{2j} \leq C_2 n (nr_n^d)^{2k-j+1} \int_{B_O(2)^{k}} \md \x \int_{(B_O(6) \setminus Q(O,\x))^{k-j+1}} \md \z \; h((O,\x)) h(\y) e^{- nr_n^d |Q(O,\x) \cup Q(\y)|} = L_{1j}+ L_{2j},
\end{equation}
where $\x = (x_2, \ldots , x_{k+1})$, $\y, \z$ are as above and
\begin{eqnarray} 
L_{1j} & = & C_2 n (nr_n^d)^{2k-j+1} \int_{B_O(2)^{k}} \md \x \int_{(B_O(6) \setminus Q(O,\x))^{k-j+1}} \md \z \; \1(|Q(O,\x) | \vee |Q(\y)| > m + \ep) \times \no \\
 & &  h((O,\x)) h(\y) e^{- nr_n^d |Q(O,\x) \cup Q(\y)|},  
\lab{eqn:defn:L_1j}
\end{eqnarray}
\begin{eqnarray} 
L_{2j} & = & C_2 n (nr_n^d)^{2k-j+1} \int_{B_O(2)^{k}} \md \x \int_{(B_O(6) \setminus Q(O,\x))^{k-j+1}} \md \z \; \1(|Q(O,\x) | \vee |Q(\y)| \leq m + \ep) \times \no \\
 & & h((O,\x)) h(\y) e^{- nr_n^d |Q(O,\x) \cup Q(\y)|}.
\lab{eqn:defn:L_2j}
\end{eqnarray}
and $\ep > 0$ is arbitrary. Using the restriction $|Q(O,\x) | \vee |Q(\y)| > m + \ep$ and substituting for $r_n$ in (\ref{eqn:defn:L_1j}) yields the bound 
\begin{equation}
L_{1j} \leq C_3 n (nr_n^d)^{2k-j+1} e^{- \f{m + \ep}{c_n} (\log n + (k-1) \log \log n + \al)}.
\lab{eqn:bound_L_1j}
\end{equation}
Choose $\eta \in (0, \f{\ep}{m}).$ Since $c_n \to m$, we can choose $\eta < \f{\ep}{m}$ sufficiently small so that $m + \ep > (1 + \eta) c_n$ for all $n$ sufficiently large. Using this in (\ref{eqn:bound_L_1j}) we obtain that as $n \to \infty$
\begin{equation}
L_{1j} \leq C_4 \f{ n (\log n)^{2k - j + 1}}{n^{1+\eta}} \to 0.
\lab{eqn:cgs_L_1j}
\end{equation}

It remains to show that $L_{2j} \to 0$. Denote by
\begin{equation}
\tilde{D}_{j,\epsilon} = \{ (\x,\z) \in B_O(2)^{k} \times (B_O(6) \setminus Q(O,\x))^{k-j+1}: |Q(O,\x) | \vee |Q(\y)| \leq m + \ep, h((O,\x)) =1, h(\y) = 1 \},
\lab{eqn:setD}
\end{equation}
the region of integration in (\ref{eqn:defn:L_2j}), where  $\x = (x_2, \ldots , x_{k+1}) \in \mR^{dk}$, $\z = (z_1,\ldots z_{k-j+1}) \in \mR^{d(k-j+1)}$, and $\y = (x_{k-j+2}, \ldots ,x_{k+1}, z_1,\ldots z_{k-j+1}) \in \mR^{d(k+1)}$. 

Since $|Q(O,\x) | \vee |Q(\y)|$ is continuous in $\tilde{D}_{j,\epsilon}$ and the minimum $m$ is achieved in $D_{j,\delta}$ for any $\delta > 0$, we can choose an $\epsilon$ small enough such that there exists a $\delta_0$ small enough with $\tilde{D}_{j,\epsilon} \subset D_{j,\delta_0}$ where $D_{j,\delta_0}$ is defined in Lemma \ref{lem:geometric}. It then follows from Lemma \ref{lem:geometric}, (\ref{eqn:defn:L_2j}), (\ref{eqn:def_beta_n}) and (\ref{def:stabilising_radii}) that for $\epsilon > 0$ sufficiently small, 
\begin{eqnarray} 
L_{2j} & \leq & C_2 n (nr_n^d)^{2k-j+1} \int_{\tilde{D}_{j,\epsilon} } \md \x \md \z \, e^{- nr_n^d (|Q(O,\x)| + |Q(\y) \setminus Q(O,\x)|)} \no \\
& \leq &  C_5 n (nr_n^d)^{2k-j+1} \int_{\tilde{D}_{j,\epsilon} } \md x \md \z \, e^{- nr_n^d (|Q(O,\x)| + \beta)} \no
\end{eqnarray}
\begin{eqnarray} 
 & \leq & C_6 \beta_n (n r_n^d)^{k-j+1} e^{- nr_n^d \beta} \leq C_7 (\log n)^{k-j+1}n^{-\beta} \to 0. 
\end{eqnarray}

This completes the proof of Theorem~\ref{thm:rips-Poisson-cgs}. \qed

 {\bf Proof of Theorem \ref{thm:rips-Poisson-cgs_finite}:}
We shall again fix $1 \leq k \leq d, p = \cR, q = \cU$ and omit these subscripts and supersctipts. Further, let $r_n = r_n^{\cR,\cU}(c_n)$. By Slutksy's lemma and Theorem \ref{thm:rips-Poisson-cgs}, it suffices to show that for any $L \geq 2$, $\EXP{J_n(r_n,L)} \to 0$.
Now, fix $L \geq 2$ and let $\Gamma \in \Delta(L,M)$ (Recall the notation in \eqref{eqn:jnkrl} and \eqref{eqn:jnkrgamma} from the proof of Theorem \ref{thm:finite_components}). We shall show that 
$$ \EXP{J_n(r_n,\Gamma)} \to 0.$$
Deriving as in the proof of Theorem \ref{thm:finite_components} (and using the same notation), we have the following bound :
$$ \EXP{J_n(r_n,\Gamma)} \leq \frac{n(n(r_n)^d)^{M-1}}{M!}\int_{(B_O(K))^{M-1}}\prod_{i=1}^Lh(\x^i)e^{-n(r_n)^d|\cup_{i=0}^LQ(\x^i)|} \md x_2 \ldots \md x_M,$$
where for all $1 \leq i \leq L$, $\x^i$ is a $(k+1)$-subset of $\{x_1,\ldots,x_M\}$ where we have set $x_1 = O$. Now, we shall break the proof into three cases. Fix $\epsilon > 0$ which will be chosen later. We shall break the integral into three cases i.e., define
\begin{eqnarray*}
A_1 &:=& \{ \x \in (B_O(K))^{M-1} : \max_{i=1,\ldots,L}|Q(\x^i)| > m_k + \epsilon \}, \\
A_2 &:= &  \{ \x \in (B_O(K))^{M-1} : \max_{i=1,\ldots,L}|Q(\x^i)| \leq m_k + \epsilon, \max_{1 \leq i < j \leq M}|x_i-x_j| > 2 \}, \\
A_3 &:=&  \{ \x \in (B_O(K))^{M-1} : \max_{i=1,\ldots,L}|Q(\x^i)| \leq m_k + \epsilon, \max_{1 \leq i < j \leq M}|x_i-x_j| \leq 2 \}.
\end{eqnarray*}
Thus, we write $ \EXP{J_n(r_n,\Gamma)} \leq I_1 + I_2 + I_3$, where $I_i$'s are defined as
$$ I_i := \frac{n(n(r_n^*)^d)^{M-1}}{M!}\int_{A_i}\prod_{i=1}^Lh(\x^i)e^{-n(r_n^*)^d|\cup_{i=0}^LQ(\x^i)|} \md \x_2 \ldots \md \x_M.$$
Now we shall show that $I_i \to 0$ as $n \to \infty$ for $i =1,2,3$ and thus complete the proof. 

First consider $I_1$. Here we have that  $\max_{i=1,\ldots,L}|Q(\x^i)| > m_k + \epsilon$.  This is the easiest of the three cases. Here using a bound similar to $L_{1j}$ in \eqref{eqn:bound_L_1j}, we can show that $I_1$ converges to $0$ as in \eqref{eqn:cgs_L_1j}.

The analysis of the remaining two cases will follow along similar lines to the bounds obtained for $L_{2j}$ in the proof of Theorem \ref{thm:rips-Poisson-cgs} using Lemma \ref{lem:geometric}. 

Next consider $I_2$. Here we have that $\max_{i=1,\ldots,L}|Q(\x^i)| \leq m_k + \epsilon$ but $\max_{1 \leq i < j \leq M}|x_i-x_j| > 2$. 

Without loss of generality, re-write $\x^1 = (O,\x) = (O,x_2, \ldots ,x_{k+1}), \x^2 = \y = (x_{k-j+2},\ldots ,x_{k+1},z_1, \ldots ,z_{k-j+1})$ for some $j \geq 1$ with $|z_1| > 2$. Setting $\z = (z_1, \ldots ,z_{k-j+1})$, we have that $(\x,\z) \in \tilde{D}_{j,\epsilon}$, where $\tilde{D}_{j,\epsilon}$ is as defined in \eqref{eqn:setD}. Now, as we argued below \eqref{eqn:setD} using continuity of $Q(.)$ as well as minimum being achieved in $D_{j,\delta}$, we again have that for $\epsilon$ small enough, there exists $\delta_0 > 0$ such that $\tilde{D}_{j,\epsilon} \subset D_{j,\delta_0}$ and hence from the geometric Lemma \ref{lem:geometric}, we have that the following inequality holds for some $\beta > 0$ :
\begin{equation}
\label{eqn:qineq}
|Q(\x^1) \cup Q(\x^2)| \geq |Q(\x^1)| + \beta.
\end{equation}
Thus we have that for some constant $C$,
$$ I_2 \leq C (n(r_n^*)^d)^{M-k-1}e^{-\beta n(r_n^*)^d} \times n (n(r_n^*)^d)^k\int_{(B_O(K))^k}e^{-n(r_n^*)^d|Q(O,x_2,\ldots,x_{k+1})|} \md x_2 \ldots \md x_{k+1}.$$
Since the convergence of the latter term on the RHS follows by Proposition \ref{prop:convergence_expectation} and the first term converges to $0$, we have that $I_2 \to 0$. 

Finally, consider $I_3$ and here we have that $\max_{i=1,\ldots,L}|Q(\x^i)| \leq m + \epsilon$ but $\max_{1 \leq i < j \leq M}|x_i-x_j| \leq 2$. 

Firstly note that this means that $M \leq d+1$ but because $\Gamma$ is a component of order at least two, $M \geq k+2$. Thus, $A_3 = \emptyset$ unless we assume that $k+2 \leq M \leq d+1$. Further, $L = \binom{M}{k+1}$ since all $(k+1)$-tuples will form $k$-simplices.

Since for all $k$-simplices $\x^i$ we have that $|Q(\x^i)| \leq m + \epsilon$ and by continuity of $|Q(\x^i)|$, there exists $\delta > 0$ (depending on $\epsilon$) such that $2 -\delta \leq |x_i-x_j| \leq 2$ for all $1 \leq i < j \leq M$. Again using geometric Lemma \ref{lem:geometric}, we have that for $\epsilon$ sufficiently small, there exists $\beta > 0$ such that \eqref{eqn:qineq} holds. Now, by proceeding as in case of $I_2$, we conclude that $I_3 \to 0$ as $n \to \infty$.   \qed 

\section*{Acknowledgements:}
The authors are thankful to Omer Bobrowski for many discussions as well as sharing drafts of his preprints. 

\bibliographystyle{plain}
\bibliography{threshold_isolated_faces_arxiv}
\end{document}